\documentclass[a4paper,11pt]{article}
\usepackage{graphicx}
\usepackage{amsfonts}
\usepackage{amsmath}
\usepackage{amssymb}
\usepackage{indentfirst}
\usepackage{titlesec}

\def\barr{\begin{array}}
\def\earr{\end{array}}
\def\bali{\begin{aligned}}
\def\eali{\end{aligned}}
\def\bearr{\begin{eqnarray}}
\def\eearr{\end{eqnarray}}

\providecommand{\play}{\displaystyle}

\providecommand{\li}{\limits}

\providecommand{\pt}{\partial}

\providecommand{\ra}{\rightarrow}

\providecommand{\da}{\downarrow}

\providecommand{\Prob}{\mathbf P}

\providecommand{\E}{\mathbf E}

\providecommand{\al}{\alpha}

\providecommand{\bt}{\beta}

\providecommand{\gm}{\gamma}

\providecommand{\dt}{\delta}

\providecommand{\Dt}{\Delta}

\providecommand{\ve}{\varepsilon}

\providecommand{\tht}{\theta}

\providecommand{\kp}{\kappa}

\providecommand{\lb}{\lambda}

\providecommand{\Lb}{\Lambda}

\providecommand{\sm}{\sigma}

\providecommand{\Om}{\Omega}

\providecommand{\R}{\mathbb R}

\providecommand{\1}{\mathbf 1}

\providecommand{\contfunc}{\mathbf{C}}

\providecommand{\grad}{\nabla}

\providecommand{\boldq}{\boldsymbol{q}}

\providecommand{\boldp}{\boldsymbol{p}}

\providecommand{\boldW}{\boldsymbol{W}}

\providecommand{\boldb}{\boldsymbol{b}}

\providecommand{\boldn}{\boldsymbol{n}}

\providecommand{\qtl}{\widetilde{q}}

\providecommand{\utl}{\widetilde{u}}

\providecommand{\vtl}{\widetilde{v}}

\providecommand{\ytl}{\widetilde{y}}

\providecommand{\lbtl}{\widetilde{\lambda}}

\providecommand{\boldqtl}{\boldsymbol{\widetilde{q}}}

\providecommand{\boldpi}{\boldsymbol{\pi}}

\providecommand{\fC}{\mathfrak{C}}

\providecommand{\fo}{\mathfrak{o}}

\begin{document}

\title{Small mass asymptotic for the motion with vanishing friction}
\author{Mark Freidlin\thanks{Department of Mathematics,
University of Maryland at College Park, mif@math.umd.edu} \ , \
Wenqing Hu\thanks{Department of Mathematics, University of Maryland
at College Park, huwenqing@math.umd.edu} \ , \ Alexander
Wentzell\thanks{Department of Mathematics, Tulane University,
wentzell@math.tulane.edu}}
\date{}

\maketitle

\begin{abstract}
We consider the small mass asymptotic (Smoluchowski-Kramers
approximation) for the Langevin equation with a variable friction
coefficient. The friction coefficient is assumed to be vanishing
within certain region. We introduce a regularization for this
problem and study the limiting motion for the 1-dimensional case and
a multidimensional model problem. The limiting motion is a Markov
process on a projected space. We specify the generator and boundary
condition of this limiting Markov process and prove the convergence.
\end{abstract}



\section{Introduction}

The Langevin equation $$\mu
\ddot{\boldq}_t^\mu=\boldb(\boldq_t^\mu)-\lb
\dot{\boldq}_t^\mu+\sm(\boldq_t^\mu)\dot{\boldW}_t \ , \
\boldq_0^\mu=\boldq\in \R^n \ , \dot{\boldq}_0^\mu=\boldp\in \R^n \
, \eqno(1.1)$$ describes the motion of a particle of mass $\mu$ in a
force field $\boldb(\boldq)$, $\boldq\in \R^n$, subjected to random
fluctuations and to a friction proportional to the velocity. Here
$\boldW_t$ is the standard Wiener process in $\R^n$, $\lb>0$ is the
friction coefficient. The vector field $\boldb(\boldq)$ and the
matrix function $\sm(\boldq)$ are assumed to be continuously
differentiable and bounded together with their first derivatives.
The matrix $a(\boldq)=(a_{ij}(\boldq))=\sm(\boldq)\sm^*(\boldq)$ is
assumed to be non-degenerate.

It is assumed usually that the friction coefficient $\lb$ is
 a positive constant. Under this assumption, one can prove that $\boldq_t^\mu$
converges in probability as $\mu \da 0$ uniformly on each finite
time interval $[0,T]$ to an $n$-dimensional diffusion process
$\boldq_t$: for any $\kp, T>0$ and any $\boldp_0^\mu=\boldp\in
\R^n$, $\boldq_0^\mu=\boldq\in \R^n$ fixed,
$$\lim\li_{\mu \da 0}\Prob\left(\max\li_{0\leq t \leq T}
|\boldq_t^\mu-\boldq_t|_{\R^d}>\kp\right)=0 \ .$$

Here $\boldq_t$ is the solution of equation

$$\dot{\boldq}_t=\dfrac{1}{\lb}\boldb(\boldq_t)+
\dfrac{1}{\lb}\sm(\boldq_t)\dot{\boldW}_t  \ , \
\boldq_0=\boldq_0^\mu=\boldq\in \R^n \ . \eqno(1.2)$$

The stochastic term in (1.2) should be understood in the It\^{o}
sense.

The approximation of $\boldq_t^\mu$ by $\boldq_t$ for $0<\mu<<1$ is
called the Smoluchowski-Kramers approximation. This is the main
justification for replacement of the second order equation (1.1) by
the first order equation (1.2). The price for such a simplification,
in particular, consists of certain non-universality of equation
(1.2): The white noise in (1.1) is an idealization of a more regular
stochastic process $\dot{\boldW}_t^\dt$ with correlation radius
$\dt<<1$ converging to $\dot{\boldW}_t$ as $\dt \da 0$. Let
$\boldq_t^{\mu,\dt}$ be the solution of equation (1.1) with
$\dot{\boldW}_t$ replaced by $\dot{\boldW}_t^\dt$. Then limit of
$\boldq_t^{\mu,\dt}$ as $\mu,\dt \da 0$ depends on the relation
between $\mu$ and $\dt$. Say, if first $\dt \da 0$ and then $\mu \da
0$, the stochastic integral in (1.2) should be understood in the
It\^{o} sense; if first $\mu \da 0$ and then $\dt \da 0$,
$\boldq_t^{\mu,\dt}$ converges to the solution of (1.2) with
stochastic integral in the Stratonovich sense. (See, for instance,
\cite{[F JSP]}.)

We considered in \cite{[FH SK PMA]} the case of a variable friction
coefficient $\lb=\lb(\boldq)$. We assumed in that work that
$\lb(\boldq)$ is smooth and $0<\lb_0\leq \lb(\boldq)\leq
\Lb<\infty$. It turns out that in this case the solution
$\boldq_t^\mu$ of (1.1) does not converge, in general, to the
solution of (1.2) with $\lb=\lb(\boldq)$, so that the
Smoluchowski-Kramers approximation should be modified. In order to
do this, we considered in \cite{[FH SK PMA]} equation (1.1) with
$\dot{\boldW_t}$ replaced by $\dot{\boldW}_t^\dt$ described above:

$$\mu \ddot{\boldq}_t^{\mu,\dt}=\boldb(\boldq_t^{\mu,\dt})-
\lb(\boldq_t^{\mu,\dt})\dot{\boldq}_t^{\mu,\dt}
+\sm(\boldq_t^{\mu,\dt})\dot{\boldW}_t^\dt \ , \
\boldq_0^{\mu,\dt}=\boldq \ , \ \dot{\boldq}_0^{\mu,\dt}=\boldp \ .
\eqno(1.3)$$

It was proved in \cite{[FH SK PMA]} that after such a
regularization, the solution of (1.3) has a limit $\boldq^\dt_t$ as
$\mu \da 0$, and $\boldq_t^\dt$ is the unique solution of the
equation obtained from (1.3) as $\mu=0$:

$$\dot{\boldq}_t^\dt=
\dfrac{1}{\lb(\boldq_t^\dt)}\boldb(\boldq_t^\dt) +
\dfrac{1}{\lb(\boldq_t^\dt)}\sm(\boldq_t^\dt)\dot{\boldW}_t^\dt
 \ , \ \boldq_0^\dt=\boldq \ . \eqno(1.4)$$

Now we can take $\dt \da 0$ in (1.4). As the result we get the
equation
$$\dot{\boldq}_t=\dfrac{1}{\lb(\boldq_t)}\boldb(\boldq_t)+
\dfrac{1}{\lb(\boldq_t)}\sm(\boldq_t)\circ \dot{\boldW}_t \ , \
\boldq_0=\boldq \ , \eqno(1.5)$$ where the stochastic term should be
understood in the Stratonovich sense. We have, for any $\dt ,\kp , T
>0$ fixed and any $\boldp_0^{\mu,\dt}=\boldp$ fixed, that
$$\lim\li_{\mu \da 0}\Prob\left(\max\li_{0\leq t \leq T}|\boldq_t^{\mu,\dt}-\boldq_t^\dt|_{\R^d}>\kp\right)=0 \ , $$
and we have $$\lim\li_{\dt \ra 0}\E \max\li_{t\in
[0,T]}|\boldq_t^\dt-\boldq_t|_{\R^d}=0  \ .
$$ So the regularization leads to a modified Smoluchowski-Kramers
equation (1.5).

In this paper we study a further generalization of the problem
considered in \cite{[FH SK PMA]}. Keeping the assumptions on uniform
boundedness and smoothness of $\lb(\bullet)$, we drop the assumption
that $0<\lb_0\leq \lb(\boldq)$ and instead assume that
$\lb(\boldq)=0$ for $\boldq\in [G] \subset \R^n$ and $\lb(\boldq)>0$
for $\boldq\in \R^n \backslash [G]$. Here $G$ is a domain in $\R^n$
and $[G]$ its closure in the standard Euclidean metric. For
simplicity of presentation we assume in the rest of this paper that
$\sm(\bullet)$ is the identity matrix. (In Section 3 we further
assume that $\boldb(\bullet)=\mathbf{0}$.) In order to use the
results of \cite{[FH SK PMA]} we introduce a further regularization
of problem (1.5). We consider the problem

$$\dot{\boldq}^\ve_t=\dfrac{1}{\lb(\boldq^\ve_t)+\ve}\boldb(\boldq^\ve_t)+
\dfrac{1}{\lb(\boldq_t^\ve)+\ve}\circ \dot{\boldW}_t \ , \
\boldq^\ve_0=\boldq \ , \ \ve>0  \eqno(1.6)$$ and we study the limit
of $\boldq_t^\ve$ as $\ve \da 0$. This limiting process can be
regarded as a limiting process of the system
$$\mu \ddot{\boldq}_t^{\mu,\dt,\ve}=\boldb(\boldq_t^{\mu,\dt,\ve})-
(\lb(\boldq_t^{\mu,\dt,\ve})+\ve)\dot{\boldq}_t^{\mu,\dt,\ve}
+\dot{\boldW}_t^\dt \ , \ \boldq_0^{\mu,\dt,\ve}=\boldq \ , \
\dot{\boldq}_0^{\mu,\dt,\ve}=\boldp \ \eqno(1.7)$$ as first $\mu \da
0$ then $\dt \da 0$ and then $\ve \da 0$.

System (1.6), in It\^{o}'s form, can be written as follows:
$$\dot{\boldq}^\ve_t=\dfrac{1}{\lb(\boldq^\ve_t)+\ve}\boldb(\boldq^\ve_t)
-\dfrac{\grad \lb(\boldq_t^\ve)}{2(\lb(\boldq_t^\ve)+\ve)^3}+
\dfrac{1}{\lb(\boldq_t^\ve)+\ve} \dot{\boldW}_t \ , \
\boldq^\ve_0=\boldq \ . \eqno(1.8)$$

However, as will be shown later, for non-compact region $[G]$, it is
sometimes more convenient to consider the projection of the above
system onto another space $\mathfrak{X}$. (In particular, in Section
3 the space $\mathfrak{X}$ is a cylinder $\mathfrak{X}=S^1\times
[a-1,b+1]$ for $a<0, b>0$.) Let us work with system (1.8) on
$\mathfrak{X}$ and compact region $[G]$. It turns out that, in the
limit, to get a Markov process with continuous trajectories, one has
to glue all the points of $[G]$ and form a projected space $\fC$.
Let the projection map be $\boldpi: \mathfrak{X} \ra \fC$. We will
prove, for the 1-dimensional case (Section 2) and a multidimensional
model problem (Section 3), that the processes
$\boldqtl_t^\ve=\boldpi(\boldq_t^\ve)$ converge weakly as $\ve \da
0$ to a continuous strong Markov process $\boldqtl_t$ on $\fC$. We
will characterize the generator of this Markov process and specify
its boundary condition. In particular, we will show that as $\ve>0$
is very small, certain mixing within $[G]$ is likely to happen for
the process $\boldq_t^\ve$. This mixing is the key mechanism that
leads to our special boundary condition. We expect that (see Section
4), within the region that the friction is vanishing, similar mixing
phenomenon will happen for the general multidimensional case.

It is worth mentioning here that some related problems are
considered in \cite{[Mochanov]}, \cite{[Mochanov-Ostrovskii]},
\cite{[Ueno I]} and \cite{[Ueno II]}. It is also interesting to note
that the limiting process for our two dimensional model problem (see
Section 3) shares some common feature with the so called Walsh's
Brownian motion (see, for example \cite{[Barlow-Pitman-Yor]}).

However, at this stage we are not able to prove, in the most general
multidimensional case (except for the 2-d model problem in Section
3), the convergence of $\boldqtl_t^\ve=\boldpi(\boldq_t^\ve)$ in
(1.8) to some Markov process $\boldqtl_t$. We will formulate a
conjecture about this in Section 4.

\section{One dimensional case}

Let us consider in this section the 1-dimensional case. Besides the
usual assumptions made in Section 1 we suppose that our friction
$\lb(\bullet)$ satisfies $\lb(q)>0$ for $q\in (-\infty, -1)\cup
(1,\infty)$. Let $\lb(q)=0$ for $q\in [-1,1]$. Equation (1.8) now
takes the following form:

$$\dot{q}_t^{\ve}=\dfrac{b(q_t^\ve)}{\lb(q_t^\ve)+\ve}-
\dfrac{\lb'(q_t^\ve)}{2(\lb(q_t^\ve)+\ve)^3}+\dfrac{1}{\lb(q_t^\ve)+\ve}\dot{W}_t
\ , \ q_0^\ve=q_0\in \R \ . \eqno(2.1)$$

We suppose that $q_0\in [a-1,b+1]$ for some $a<0<b$. The process
$q_t^\ve$ is supposed to be stopped once it hits $q=a-1$ or $q=b+1$.

Our goal is to study the asymptotic behavior of (2.1) as $\ve \da
0$. To this end we shall write the process (2.1) as a strong Markov
process subject to a generalized second order differential operator
in the form $D_{v^\ve}D_{u^\ve}$ (see \cite{[Feller]},
\cite{[Dynkin]}, \cite{[Mandl]}). We have

$$u^\ve(q)=\int_0^q(\lb(x)+\ve)\exp\left(-2\int_0^x b(y)(\lb(y)+\ve)dy\right)dx  \ , \eqno(2.2)$$

$$v^\ve(q)=2\int_0^q(\lb(x)+\ve)\exp\left(2\int_0^x b(y)(\lb(y)+\ve)dy\right)dx  \ . \eqno(2.3)$$

For fixed $\ve>0$, the functions $u^\ve$ and $v^\ve$ are strictly
increasing functions in their arguments. As $\ve\da 0$, they will
converge uniformly on finite intervals to the functions $u$ and $v$
defined by

$$u(q)=\int_0^q\lb(x)\exp\left(-2\int_0^x b(y)\lb(y)dy\right)dx  \ , \eqno(2.4)$$

$$v(q)=2\int_0^q\lb(x)\exp\left(2\int_0^x b(y)\lb(y)dy\right)dx  \ . \eqno(2.5)$$

The functions $u$ and $v$ are strictly increasing outside the
interval $[-1,1]$ and have constant stretches on $[-1,1]$.

\

Consider a projection map $\pi$: we let $\pi([-1,1])=0$ and
$\pi(q)=q+1$ for $q<-1$ and $\pi(q)=q-1$ for $q>1$. Consider the
process $\qtl^\ve_t=\pi(q_t^\ve)$. Process $\qtl^\ve_t$ for fixed
$\ve>0$, in general, is \textit{not} a Markov process.

\

Let us define two functions $\utl$ and $\vtl$ as follows:
$\utl(\qtl)=u(\qtl-1)$ for $\qtl<0$ and $\utl(\qtl)=u(\qtl+1)$ for
$\qtl>0$ and $\utl(0)=u(1)=u(-1)=0$; $\vtl(\qtl)=v(\qtl-1)$ for
$\qtl<0$ and $\vtl(\qtl)=v(\qtl+1)$ for $\qtl>0$ and
$\vtl(0)=v(1)=v(-1)=0$. Here the functions $u$ and $v$ are defined
in (2.4), (2.5). The functions $\utl$ and $\vtl$ are continuous
strictly increasing functions on $[a,b]$.

Define a Markov process $\qtl_t$ on $[a,b]$ as follows. The
generator $A$ of $\qtl_t$ is $A=D_{\vtl}D_{\utl}$. The domain of
definition $D(A)$ of operator $A$ consists of all functions $f$ that
are continuous on $[a,b]$, are twice continuously differentiable in
$\qtl \in [a,b]\backslash\{0\}$, with finite limit $\lim\li_{\qtl\ra
0}Af(\qtl)$ (taken as the value of $Af(0)$) and finite one-sided
limits $\lim\li_{\dt \da
0}\dfrac{f(\dt)-f(0)}{\utl(\dt)-\utl(0)}\equiv
D_{\utl}^+f(0)=D_{\utl}^-f(0)\equiv \lim\li_{\dt \da
0}\dfrac{f(0)-f(-\dt)}{\utl(0)-\utl(-\dt)}$. Also we have
$\lim\li_{\qtl\ra a}Af(\qtl)=\lim\li_{\qtl \ra b}Af(\qtl)=0$ (taken
as the value of $Af(a)$ and $Af(b)$).

\

\textbf{Lemma 2.1.} \textit{There exists the Markov process $\qtl_t$
on $[a,b]$.}

\

\textbf{Proof.} The existence of such a process could be checked
similarly as in \cite[Section 2]{[FW fish paper]}. For the sake of
completeness and comparison with results in the next section we
shall check it here. To this end we use an equivalent formulation of
the Hille-Yosida theorem (see \cite[Section 2]{[FW fish paper]} also
\cite[Theorem 2]{[Wentzell boundary condition]}). We check three
conditions.

$\bullet$ The domain $D(A)$ is dense in the space
$\mathbf{C}([a,b])$. This is because we can approximate every
continuous function $f$ with one that is constant in a neighborhood
of $0$. After that in the interior part of the intervals $[a,0)$ and
$(0,b]$, at a positive distance from 0, with a smooth function. The
approximating smooth function satisfy our boundary conditions since
$Af(0)=D_{\utl}^+f(0)=D_{\utl}^-f(0)=0$.

$\bullet$ The maximum principle: if $f\in D(A)$ and the function $f$
reaches its maximum at a point $x_0\in [a,b]$, then $Af(x_0)\leq 0$.
If $x_0\neq 0$ we have $f'(x_0)=0$ and $f''(x_0)\leq 0$ and
$$D_{\vtl}D_{\utl}f(x_0)=\dfrac{f''(x_0)}{\vtl'(x_0)\utl'(x_0)}
-\dfrac{\utl''(x_0)}{\vtl'(x_0)(\utl'(x_0))^2}f'(x_0)\leq 0 \ .$$ If
the maximum is achieved at $0$, we consider the expansion

$$f(x)=f(0)+D_{\utl}f(0)(\utl(x)-\utl(0))
+(Af(0)+o(1))\int_0^x(\vtl(y)-\vtl(0))d\utl(y) \ .$$

The last integral is $O(\utl(x)\vtl(x))$ as $x\ra 0$. Since
$D_{\utl}^-f(0)\geq 0$ and $D_{\utl}^+f(0)\leq 0$, by our boundary
conditions at $0$ we get $D_{\utl}f(0)=0$. This implies that
$Af(0)\leq 0$.

$\bullet$ Existence of solution $f\in D(A)$ of $\lb f-A f=F$ for all
$F\in \mathbf{C}([a,b])$. On each of the intervals $[a,0)$ and
$(0,b]$ the general solution of equation $\lb f - D_{\vtl}D_{\utl} f
= F$, $F\in \mathbf{C}([a,b])$ can be written as

$$f^{\pm}(q)=\widehat{f}^{\pm}(q)+G^{\pm}(q) \ .$$

Here $\widehat{f}^{\pm}(q)$ satisfy the equation $\lb
\widehat{f}^{\pm}-D_{\vtl}D_{\utl} \widehat{f}^{\pm}=F$,
$\widehat{f}^+(0+)=0$ (or $\widehat{f}^-(0-)=0$),
$D_{\utl}^+\widehat{f}^+(0)=0$ (or $D_{\utl}^-\widehat{f}^-(0)=0$)
and $G^\pm(q)$ satisfy the equation $\lb G^\pm-D_{\vtl}D_{\utl}
G^\pm=0$, $G^+(0+)=k_1^+$ (or $G^-(0-)=k_1^-$),
$D_{\utl}^+G^+(0)=k_2^+$ (or $D_{\utl}^-G^-(0)=k_2^-$). Here
$k_1^{\pm}$ and $k_2^{\pm}$ are constants. Our boundary condition
gives $k_1^{+}=k_1^{-}$ and $k_2^{+}=k_2^{-}$. The boundary
condition $D_{\utl}D_{\vtl}f^{+}(a)=D_{\utl}D_{\vtl}f^{-}(b)=0$
singles out a unique $f\in D(A)$. $\square$

\

We have

\

\textbf{Theorem 2.1.} \textit{As $\ve \da 0$, for fixed $T>0$, the
process $\qtl_t^\ve$ converges weakly in the space
$\mathbf{C}_{[0,T]}([a,b])$ to the process $\qtl_t$.}

\

The proof of this Theorem is based on an application of the
machinery developed in \cite[Ch.8]{[FW book]}, \cite{[FW AMS]} and
\cite{[FW fish paper]}. We shall use the following lemma, which is
the Lemma 3.1 of \cite[Ch.8, page 301]{[FW book]}. We formulate it
here in the terminology that meets our purpose.

\

\textbf{Lemma 2.2.} \textit{Let $M$ be a metric space; $Y$, a
continuous mapping $M \mapsto Y(M)$, $Y(M)$ being a complete
separable metric space. Let $(X_t^\ve, \Prob_x^\ve)$ be a family of
Markov processes in $M$; suppose that the process $Y(X_t^\ve)$ has
continuous trajectories. Let $(y_t, \Prob_y)$ be a Markov process
with continuous paths in $Y(M)$ whose infinitesimal operator is $A$
with domain of definition $D(A)$. Let $T>0$. Let us suppose that the
space $\mathbf{C}_{[0,T]}(Y(M))$ of continuous functions on $[0,T]$
with values in $Y(M)$ is taken as the sample space, so that the
distribution of the process in the space of continuous functions is
simply $\Prob_y$. Let $\Psi$ be a subset of the space
$\mathbf{C}_{[0,\infty)}(Y(M))$ such that for measures $\mu_1$,
$\mu_2$ on $Y(M)$ the equality $\play{\int F d\mu_1=\int F d\mu_2}$
for all $F\in \Psi$ implies $\mu_1=\mu_2$. Let $D$ be the subset of
$D(A)$ such that for every $F\in \Psi$ and $\lb>0$ the equation $\lb
f - A f=F$ has a solution $f\in D$.}

\textit{Suppose that for every $x\in M$ the family of distributions
$\mathbf{Q}_x^\ve$ of $Y(X_{\bullet}^\ve)$ in the space
$\mathbf{C}_{[0,T]}(Y(M))$ corresponding to the probabilities of
$\Prob_x^\ve$ is weakly pre-compact; and that for every compact $K
\subset Y(M)$, for every $f\in D$ and every $\lb > 0$,}

$$\E_x^\ve\int_0^\infty e^{-\lb t}[\lb f(Y(X_t^\ve))-A f(Y(X_t^\ve))]dt\ra f(Y(x))$$
\textit{as $\ve \da 0$ uniformly in $x\in Y^{-1}(K)$.}

\textit{Then $\mathbf{Q}_x^\ve$ converges weakly as $\ve \da 0$ to
the probability measure $\Prob_{Y(x)}$.}

\

\textbf{Proof of Theorem 2.1.} Making use of Lemma 2.2, we take the
metric space $M=[a-1,b+1]$ and the mapping $Y=\pi$. The space
$Y(M)=\pi([a-1,b+1])=[a,b]$. We take the process $q_t^\ve$ as
$(X_t^\ve, \Prob_x^\ve)$. We take the process $\qtl_t$ as $(y_t,
\Prob_y)$.

Let $\Psi$ be the space of all continuous bounded functions in
$[a,b]$ which are once continuously differentiable inside $[a,0)$
and $(0,b]$, with bounded derivatives. The space $D\subset D(A)$
consists of those functions $f\in D(A)$ such that they are
continuous and bounded in $[a,b]$ and are three times continuously
differentiable inside $[a,0)$ and $(0,b]$, with bounded derivatives
up to the third order.

Pre-compactness of the family of distributions of the process
$\{\qtl^\ve_{\bullet}\}_{\ve>0}$ is checked in Lemma 2.4. What
remains to do is to check that for every compact $K \subset [a,b]$,
for every $f\in D$ and every $\lb > 0$,
$$\E_{q_0}\left[ \int_0^\infty e^{-\lb t}[\lb f(\pi(q_t^\ve))-A f(\pi(q_t^\ve))]dt- f(\pi(q_0))\right]\ra 0$$
as $\ve \da 0$ uniformly in $q_0\in \pi^{-1}(K)$. This is done in
Lemma 2.5. This finishes the proof of Theorem 2.1. $\square$

\

For positive $\dt$ small enough, let $G(\dt)=[a-1,
-1-\dt]\cup[1+\dt,b+1]$. Let $0<\dt'<\dt$. Let $C(\dt')=\{-1-\dt',
1+\dt'\}$. We introduce a sequence of stopping times $\tau_0\leq
\sm_0<\tau_1<\sm_1<\tau_2<\sm_2<...$ by

$$\tau_0=0 \ ,
\ \sm_n=\min\{t\geq \tau_n, q_t^\ve\in G(\dt)\} \ , \
\tau_n=\min\{t>\sm_{n-1}: q_t^\ve\in C(\dt')\} \ .$$

This is well-defined up to some $\sm_k$ ($k \geq 0$) such that
$$\Prob_{q_{\sm_k}^\ve}(q_{t+\sm_k}^\ve \text{ hits } a-1 \text{ or } b+1
\text{ before it hits } -1-\dt' \text{ or } 1+\dt')=1 \ .$$

We will then define $\tau_{k+1}=\min\{t>\sm_k: q_t^\ve=a-1 \text{ or
} b+1\}$. And we define
$\tau_{k+1}<\sm_{k+1}=\tau_{k+1}+1<\tau_{k+2}=\tau_{k+1}+2<\sm_{k+2}=\tau_{k+1}+3<...$
and so on.

We have $\lim\li_{n \ra \infty}\tau_n=\lim\li_{n \ra
\infty}\sm_n=\infty$. And we have obvious relations
$q_{\tau_n}^\ve\in C(\dt')$, $q_{\sm_n}^\ve\in C(\dt)$ for $1\leq n
\leq k$ (as long as $k \geq 1$, if $k=0$ the process may start from
$G(\dt)$ and goes directly to $a-1$ or $b+1$ without touching
$C(\dt')$ and is stopped there, or it may start from
$(-1-\dt,1+\dt)$, reaches $\{-1-\dt, 1+\dt\}$ first and then goes
directly to $a-1$ or $b+1$ without touching $C(\dt')$ and is stopped
there). Also, for $n\geq k+1$ we have
$q_{\tau_n}^\ve=q_{\sm_n}^\ve=a-1 \text{ or } b+1$. If
$q_0^\ve=q_0\in G(\dt)$, then we have $\sm_0=0$ and $\tau_1$ is the
first time at which the process $q_t^\ve$ reaches $C(\dt')$ or
$\{a-1,b+1\}$.

\

Now we check weak pre-compactness of the family of distributions of
the processes $\{\qtl^\ve_t\}_{\ve>0}$. To this end we need the
following lemma, which is Lemma 5.1 in \cite{[FW fish paper]}. We
formulate it using our terminology.

\

\textbf{Lemma 2.3.} \textit{Let $\qtl^{\ve,\dt}_{\bullet}$ for every
$\ve>0$, $\dt>0$, be a random element in}
$\textbf{C}_{[0,T]}([a,b])$ \textit{such that $\max\li_{0\leq t \leq
T}|\qtl^\ve_t-\qtl^{\ve,\dt}_t|\leq \dt$ on the whole probability
space. If for every positive $\dt$ the family of distributions of
$\qtl^{\ve,\dt}_{\bullet}$, $\ve > 0$, is tight, then the family of
distributions of $\qtl^\ve_{\bullet}$ is pre-compact.}

\

Now we have

\

\textbf{Lemma 2.4.} \textit{The family of distributions of
$\{\qtl^\ve_{\bullet}\}_{\ve > 0}$ is pre-compact.}

\

\textbf{Proof.} Let $\dt'=\dt/2$ so that we need only one parameter
$\dt$. Between the times $\sm_{i-1}$ and $\tau_i$ the process
$q_t^\ve$ is either in $[a,-1-\dt/2)$ or in $(1+\dt/2,b]$, and for
$\sm_{i-1}\leq t<t'< \tau_i$ we have
$|\qtl_t^\ve-\qtl_{t'}^\ve|=|q_t^\ve-q_{t'}^\ve|$. Since we have

$$q_{t}^\ve-q_{t'}^\ve=\int_t^{t'}\left[\dfrac{b(q_s^\ve)}{\lb(q_s^\ve)+\ve}
-\dfrac{\lb'(q_s^\ve)}{2(\lb(q_t^\ve)+\ve)^3}\right]ds+\int_t^{t'}\dfrac{1}{\lb(q_s^\ve)+\ve}dW_s
\ ,$$ we can estimate $$\E|q_t^\ve-q_{t'}^\ve|^4\leq K(\dt)|t-t'|^2
\ .$$ The constant $K(\dt)$ is independent of $\ve$ provided that
$\ve$ is small. Now we let

$$Z_t^{\ve,\dt}=\int_0^t \1_{G(\dt/2)}(q_s^\ve)\left[\dfrac{b(q_s^\ve)}{\lb(q_s^\ve)+\ve}
-\dfrac{\lb'(q_s^\ve)}{2(\lb(q_t^\ve)+\ve)^3}\right]ds+
\int_0^t\1_{G(\dt/2)}(q_s^\ve)\dfrac{1}{\lb(q_s^\ve)+\ve}dW_s \ .$$

From the above estimate we see that $Z_t^{\ve,\dt}$ for fixed
$\dt>0$ is tight. The trajectories of these stochastic processes
satisfy the H\"{o}lder condition
$|Z_t^{\ve,\dt}-Z_{t'}^{\ve,\dt}|\leq H^{\ve,\dt}|t-t'|^{1/5}$ where
$H^{\ve,\dt}$ are random variables with $\E(H^{\ve,\dt})^4$ bounded
by the same $K(\dt)$.

For $i \geq 1$ if $q_{\tau_i}^\ve\in C(\dt/2)$ and $q_{\sm_i}^\ve\in
C(\dt)$ then between the times $\tau_i$ and $\sm_i$ ($\leq T$) the
process $q_t^\ve$ travels a distance at least $\dt/2$ and at least
this distance in $G(\dt/2)$ on the same interval either
$[a,-1-\dt/2)$ or $(1+\dt/2,b]$. By our estimate on H\"{o}lder
continuity of $Z_t^{\ve,\dt}$ this implies that $\sm_i-\tau_i\geq
\left(\dfrac{\dt}{4H^{\ve,\dt}}\right)^5$, $i \geq 1$. If
$q_{\tau_i}^\ve \in \{a-1,b+1\}$ then by our definition of the
stopping time $\sm_i=\tau_i+1$ we can choose $\dt$ small enough such
that the above inequality also holds.

Now we shall define the process $\qtl_t^{\ve,\dt}$ as follows.

$\bullet$ For $\sm_{i-1}\leq t \leq \tau_i$ we take
$\qtl_t^{\ve,\dt}=\qtl_t^\ve$.

$\bullet$ For $\tau_0\leq t \leq \sm_0$ we take
$\qtl_t^{\ve,\dt}=\qtl_{\sm_0}^\ve$. This gives $\max\li_{\tau_0\leq
t \leq \sm_0}|\qtl_t^{\ve,\dt}-\qtl_t^\ve|=\max\li_{\tau_0\leq t
\leq \sm_0}|\qtl_{\sm_0}^\ve-\qtl_t^\ve|\leq \dt$.

$\bullet$ If $\tau_i<T<\sm_i$ we take
$\qtl_t^{\ve,\dt}=\qtl_{\tau_i}^\ve$ for $\tau_i \leq t \leq T$.
This gives $\max\li_{\tau_i\leq t \leq
T}|\qtl_t^{\ve,\dt}-\qtl_t^\ve|=\max\li_{\tau_i\leq t \leq
T}|\qtl_{\tau_i}^\ve-\qtl_t^\ve|\leq \dt/2$.

$\bullet$ If $\sm_i\leq T$. In this case if $\qtl_{\tau_i}^\ve$ and
$\qtl_{\sm_i}^\ve$ are within a distance $\leq \dt$ from $0$, we
define $\qtl_{\frac{\tau_i+\sm_i}{2}}^{\ve,\dt}=0$,
$$\qtl_{t}^{\ve,\dt}= \left(1-\dfrac{2(t-\tau_i)}{\sm_i-\tau_i}\right)
\qtl_{\tau_i}^\ve \text{ for } \tau_i\leq t \leq
\dfrac{\tau_i+\sm_i}{2} \ ,$$ $$\qtl_{t}^{\ve,\dt}=
-\left(1-\dfrac{2(t-\tau_i)}{\sm_i-\tau_i}\right) \qtl_{\sm_i}^\ve
\text{ for }  \dfrac{\tau_i+\sm_i}{2}\leq t \leq \sm_i \ .$$

Since this is just a linear interpolation it is clear that in this
case we have $\max\li_{\tau_i\leq t \leq
\sm_i}|\qtl^{\ve,\dt}_t-\qtl^\ve_t|\leq 2\dt$. Within this time
interval $\tau_i \leq t<t' \leq \sm_i$, $i\geq 1$ we have
$$|\qtl_t^{\ve,\dt}-\qtl_{t'}^{\ve,\dt}|\leq\dfrac{\dt}{|\sm_i-\tau_i|}|t-t'|\leq  \dfrac{\dt}{(\min\li_{i\geq 1}
 \dfrac{1}{2}|\sm_i-\tau_i|)^{1/5}}|t-t'|^{1/5}
\leq 2^{11/5}H^{\ve,\dt}|t-t'|^{1/5} \ .$$

Another possibility is that $q_{\sm_i}^\ve=q_{\tau_i}^\ve= a-1
\text{ or } b+1$. In this case we define
$\qtl_t^{\ve,\dt}=\qtl_t^\ve$ for $\tau_i\leq t <\sm_i$.

On the whole interval $0\leq t<t' \leq T$ we have
$|\qtl_t^{\ve,\dt}-\qtl_{t'}^{\ve,\dt}|\leq
(2^{11/5}+2)H^{\ve,\dt}|t'-t|^{1/5}$ for $|t'-t|\leq
\left(\dfrac{\dt}{4H^{\ve,\dt}}\right)^5$. This means that for fixed
$\dt>0$ we have the tightness of the family of distributions of
$\qtl_t^{\ve,\dt}$ in the space $\mathbf{C}_{[0,T]}([a,b])$. Since
we checked $\max\li_{0\leq t \leq
T}|\qtl_t^{\ve,\dt}-\qtl_t^\ve|\leq 2\dt$, by using Lemma 2.3 with
$2\dt$ instead of $\dt$ we get the pre-compactness of the family of
distributions of $\qtl_t^\ve$ in $\mathbf{C}_{[0,T]}([a,b])$.
$\square$

\

The proof of the next Lemma 2.5 is based on Lemmas 2.6-2.10. Within
the proof of this lemma and the auxiliary Lemmas 2.6-2.10, we will
take $\ve\da 0$, $\dt=\dt(\ve)\da 0$, $\dt'=\dt'(\ve)\da 0$ in an
asymptotic order such that $0<\ve<<\dt'<<\dt$. Although not very
precise, but for simplicity of presentation we will just refer this
choice of order as first $\ve \da 0$, then $\dt' \da 0$ and then
$\dt \da 0$. It could be checked that such an order of taking limit
does not alter the validity of the result.

Throughout the rest of this section and next section when we use
symbols $U$, $V$, $M_i$, $C_i$, $A_i$, etc., they are referring to
some positive constants. We will not point out this explicitly
unless some special properties of the implied constants are
stressed. Also we sometimes use the same letter for constants in
different estimates.

\

\textbf{Lemma 2.5.} \textit{For every compact $K \subset [a,b]$, for
every $f\in D$ and every $\lb > 0$,}

$$\E_{q_0}\left[ \int_0^\infty e^{-\lb t}[\lb f(\pi(q_t^\ve))-A f(\pi(q_t^\ve))]dt- f(\pi(q_0))\right]\ra 0$$
\textit{as $\ve \da 0$ uniformly in $q_0\in \pi^{-1}(K)$.}

\

\textbf{Proof.} The above expectation can be written as

$$\begin{array}{l}
\play{\E_{q_0}\left[\sum\li_{n=0}^\infty
\left[\int_{\tau_n}^{\sm_n}e^{-\lb t}[\lb
f(\pi(q_t^\ve))-Af(\pi(q_t^\ve))]dt+e^{-\lb
\sm_n}f(\pi(q_{\sm_n}^\ve))-e^{-\lb\tau_n}f(\pi(q_{\tau_n}^\ve))\right]+\right.}
\\
\play{\ \ \  \left. \sum\li_{n=0}^\infty
\left[\int_{\sm_n}^{\tau_{n+1}}e^{-\lb t}[\lb
f(\pi(q_t^\ve))-Af(\pi(q_t^\ve))]dt+e^{-\lb
\tau_{n+1}}f(\pi(q_{\tau_{n+1}}^\ve))-e^{-\lb\sm_n}f(\pi(q_{\sm_n}^\ve))\right]\right]}
\\
\play{=\E_{q_0} \left[\sum\li_{n=0}^\infty e^{-\lb
\tau_n}\psi_1^\ve(q_{\tau_n}^\ve)+\sum\li_{n=0}^\infty e^{-\lb
\sm_n} \psi_2^\ve(q_{\sm_n}^\ve)\right]} \ ,
\end{array} \eqno(2.6)$$
where

$$\psi_1^\ve(q)=\E_q\left[\int_{0}^{\sm_0}e^{-\lb t}[\lb
f(\pi(q_t^\ve))-Af(\pi(q_t^\ve))]dt+e^{-\lb
\sm_0}f(\pi(q_{\sm_0}^\ve))\right]-f(\pi(q)) \ , \eqno(2.7)$$

$$\psi_2^\ve(q)=\E_q\left[\int_{\sm_0}^{\tau_1}e^{-\lb t}[\lb
f(\pi(q_t^\ve))-Af(\pi(q_t^\ve))]dt+e^{-\lb
\tau_1}f(\pi(q_{\tau_1}^\ve))\right]-f(\pi(q)) \ . \eqno(2.8)$$

We used the strong Markov property of $q_t^\ve$. Since for $n \geq
k+1$ we have
$\psi_1^\ve(q^\ve_{\tau_n})=\psi_2^\ve(q^\ve_{\sm_n})=0$ we can
assume that the function $\psi_2^\ve$ is taken at a point on
$G(\dt)\backslash\{a-1,b+1\}$ and the expectation is determined by
the values of the process $q_t^\ve$ in one of the intervals either
$(1+\dt',b+1]$ or $[a-1, -1-\dt')$. We will prove, in Lemma 2.6,
that under our specified asymptotic order we can have
$|\psi_2^\ve(q)|\leq (\utl(\dt)-\utl(-\dt))^2$ as $\ve \da 0$.

We can assume that the function $\psi_1^\ve$ is taken at a point in
$[-1-\dt',1+\dt']$ (in the case when $n=0$ and $q_0^\ve\in G(\dt)$,
we also have $\psi_1^\ve(q_0)=0$). We can write

$$\begin{array}{l}
\psi_1^\ve(q) \\
\play{=\left(\E_q
f(\pi(q_{\sm_0}^\ve))-f(\pi(q))\right)-\E_q(1-e^{-\lb
\sm_0})f(\pi(q_{\sm_0}^\ve))+ \E_q\int_{0}^{\sm_0}e^{-\lb t}[\lb
f(\pi(q_t^\ve))-Af(\pi(q_t^\ve)))]dt}
\\
\play{=(I)^\ve(q)+(II)^\ve(q)+(III)^\ve(q) \ .}
\end{array} \eqno(2.9)$$

We are going to prove, in Lemma 2.8, that for
$q\in[-1-\dt',1+\dt']$, for a function $f\in D$ we can have the
estimate $|(I)^\ve(q)|\leq M_1(\utl(\dt)-\utl(-\dt))^2$.

In Lemma 2.9 we will show that $\E_q\sm_0\leq
M_1(\utl(\dt)-\utl(-\dt))(\vtl(\dt)-\vtl(-\dt))$ and $\E_q(1-e^{-\lb
\sm_0})\leq M_1(\utl(\dt)-\utl(-\dt))(\vtl(\dt)-\vtl(-\dt))$ so that
$|(II)^\ve(q)|+|(III)^\ve(q)|<M_1(\utl(\dt)-\utl(-\dt))(\vtl(\dt)-\vtl(-\dt))$
for $q\in [-1-\dt',1+\dt']$.

These estimates show that
$$|\psi^\ve_1(q)|<(\utl(\dt)-\utl(-\dt))^2+M_1(\utl(\dt)-\utl(-\dt))(\vtl(\dt)-\vtl(-\dt))$$
for all $q\in [-1-\dt',1+\dt']$.

As we only consider the arguments $q_{\tau_n}^\ve$ of $\psi_1^\ve$
in (2.6) being in $[-1-\dt',1+\dt']$ starting with $n=1$ (otherwise
$\psi_1^\ve=0$), we have, by strong Markov property of $q_t^\ve$,
that

$$\begin{array}{l}
\play{\left|\E_{q_0}\sum\li_{n=1}^\infty e^{-\lb
\tau_n}\psi_1^\ve(q_{\tau_n}^\ve)\right|}\
\\
\play{ \leq
(\utl(\dt)-\utl(-\dt))^2+M_1(\utl(\dt)-\utl(-\dt))(\vtl(\dt)-\vtl(-\dt))\sum\li_{n=1}^\infty
\E_{q_0}e^{-\lb \tau_n}}
\\
\play{ \leq
(\utl(\dt)-\utl(-\dt))^2+M_1(\utl(\dt)-\utl(-\dt))(\vtl(\dt)-\vtl(-\dt))\sum\li_{n=1}^\infty
\left(\sup\li_{q\in G(\dt)}\E_{q}e^{-\lb \tau_1}\right)^{n-1}} \ .
\end{array}$$

We will show, in Lemma 2.10, that $\E_{q}e^{-\lb
\tau_1}<1-M_2\utl(\dt)\wedge (-\utl(-\dt))$ for all $q\in G(\dt)$.
Since as $\dt \da 0$ we have $0<M_2\leq
\left|\dfrac{\utl(\dt)}{-\utl(-\dt)}\right|\leq M_3<\infty$, we have
$$\begin{array}{l}
\play{\left|\E_{q_0}\sum\li_{n=1}^\infty e^{-\lb
\tau_n}\psi_1^\ve(q_{\tau_n}^\ve)\right|}
\\
\play{\leq
((\utl(\dt)-\utl(-\dt))^2+M_1(\utl(\dt)-\utl(-\dt))(\vtl(\dt)-\vtl(-\dt)))\dfrac{1}{M_2(\utl(\dt))\wedge
(-\utl(-\dt))}}\ra 0  \end{array}$$ as $\dt \da 0$. For $n=0$ the
expectation $\E_{q_0}\psi^\ve_1(q_0^\ve)$ is small as $\ve$ is
small.

For the second term in (2.6) we can estimate

$$\begin{array}{l}
\play{\left|\sum\li_{n=0}^\infty \E_{q} e^{-\lb \sm_n}
\psi_2^\ve(q^\ve_{\sm_n})\right| \leq \sum\li_{n=0}^\infty \E_{q}
e^{-\lb\sm_n}|\psi_2^{\ve}(q)| \leq \sum\li_{n=0}^\infty \E_{q}
e^{-\lb\tau_n}|\psi_2^{\ve}(q)|}
\\
\play{\leq (1+\dfrac{M_4}{(\utl(\dt))\wedge
(-\utl(-\dt))})(\utl(\dt)-\utl(-\dt))^2}
\end{array}
$$
which converges to $0$ as $\ve \da 0$. This proves this lemma.
$\square$

\

\textbf{Lemma 2.6.} \textit{We have, for $q\in G(\dt)$, as $\ve$ is
small, that $|\psi_2^\ve(q)|\leq (\utl(\dt)-\utl(-\dt))^2$.}

\

\textbf{Proof.} For the initial point $q\in G(\dt)$ and the time
interval $0\leq t \leq \tau_1$ the trajectory of $q_t^\ve$ is
traveling in one of the intervals either $[1+\dt',1+b]$ or
$[a-1,-1-\dt']$. Without loss of generality let us assume that $q\in
[1+\dt,1+b]$ and we are traveling in the interval $[1+\dt',1+b]$.
Let $\qtl=\pi(q)$. Let $B(\qtl)=b(\qtl+1)$ and
$\Lb(\qtl)=\lb(\qtl+1)$. Let us extend the function $\Lb(\bullet)$
to the whole line $\R$. The extended function
$\widehat{\Lb}(\bullet)$ is smooth, bounded, with uniformly bounded
derivatives and such that $\widehat{\Lb}(x)\geq \min\li_{q\in
[1+\dt',1+b]}\lb(q)$, $\widehat{\Lb}(x)=\lb(1+x)$ for $x\in
[\dt',b]$.

Let the process $\widehat{\qtl}^\ve_t$ be subject to the stochastic
differential equation

$$\dot{\widehat{\qtl}}^\ve_t=\dfrac{B(\widehat{\qtl}_t^\ve)}{\widehat{\Lb}(\widehat{\qtl}_t^\ve)+\ve}
-\dfrac{\widehat{\Lb}'(\widehat{\qtl}_t^\ve)}{2(\widehat{\Lb}(\widehat{\qtl}_t^\ve)+\ve)^3}
+\dfrac{1}{\widehat{\Lb}(\widehat{\qtl}_t^\ve)+\ve}\dot{W}_t \ , \
\widehat{\qtl}_0^\ve=\qtl \ , 0\leq t <\infty \ .$$

We introduce a stochastic process $\widehat{\qtl}_t$,
$\widehat{\qtl}_0=\qtl$ with generator $\widehat{A}$, subject to the
stochastic differential equation

$$\dot{\widehat{\qtl}}_t=\dfrac{B(\widehat{\qtl}_t)}{\widehat{\Lb}(\widehat{\qtl}_t)}
-\dfrac{\widehat{\Lb}'(\widehat{\qtl}_t)}{2\widehat{\Lb}^3(\widehat{\qtl}_t)}
+\dfrac{1}{\widehat{\Lb}(\widehat{\qtl}_t)}\dot{W}_t \ , \
\widehat{\qtl}_0=\qtl \ , 0\leq t <\infty \ .$$

Notice that the modified generator $\widehat{A}$ agrees with $A$
before the process $\qtl_t^\ve$ reaches $\qtl_{\tau_1}^\ve$. And
before the time $\tau_1$ the process $\widehat{\qtl}_t^\ve$ agrees
with the process $\qtl_t^\ve$. Therefore we have,

$$\psi_2^\ve(q)=\E_{\qtl}\left[\int_0^{\tau_1}e^{-\lb t}[\lb f(\widehat{\qtl}_t^\ve)
-\widehat{A}f(\widehat{\qtl}_t^\ve)]dt-e^{-\lb
\tau_1}f(\widehat{\qtl}_{\tau_1}^\ve)\right]-f(\qtl) \ .$$

It is clear by It\^{o}'s formula that we have (also see,
\cite[Section 2]{[FW Diffusion process on a graph]}), for the
stopping time $\tau_1$,

$$\E_{\qtl}\left[\int_0^{\tau_1}e^{-\lb t}[\lb
f(\widehat{\qtl}_t)-\widehat{A}f(\widehat{\qtl}_t)]dt-e^{-\lb
\tau_1}f(\widehat{\qtl}_{\tau_1})\right] -f(\qtl)=0 \ .$$

Notice that the function $f\in D\subset D(A)$ is three times
continuously differentiable in $[\dt',b]$. This gives the estimate
that for some positive $U,V>0$ and $T=T(\ve)$ we have

$$\begin{array}{l}
|\psi_2^\ve(q)|
\\
\play{=\left|\E_{\qtl}\int_0^{\tau_1} e^{-\lb t}[\lb
(f(\widehat{\qtl}_t^\ve)-f(\widehat{\qtl}_t))-
(\widehat{A}f(\widehat{\qtl}_t^\ve)-\widehat{A}f(\widehat{\qtl}_t))]dt
-e^{-\lb
\tau_1}(f(\widehat{\qtl}_{\tau_1}^\ve)-f(\widehat{\qtl}_{\tau_1}))\right|}
\\
\play{\leq \E_{\qtl}\left(\int_0^{T(\ve)} \lb e^{-\lb t}dt
\left(\text{Lip}(f)\right)\cdot|\widehat{\qtl}_t^\ve-\widehat{\qtl}_t|
+ \right.}
\\
\play{\left. \ \ \ \ \ \ \ \ \ \ \int_0^{T(\ve)}  e^{-\lb t}dt
\left(\text{Lip}(Af)\right)\cdot|\widehat{\qtl}_t^\ve-\widehat{\qtl}_t|
+
\left(\text{Lip}(f)\right)\cdot|\widehat{\qtl}_{\tau_1}^\ve-\widehat{\qtl}_{\tau_1}|\1\left(\tau_1\leq
T(\ve)\right) \right) +}
\\
\play{\ \ \ \ \ \ \ \ \ \ \ \ \ \ \ \ \ \ \ V \Prob(\tau_1\geq
T(\ve))}
\\
\play{\leq U \left(\max\li_{0\leq t \leq
T(\ve)}\E_{\qtl}|\widehat{\qtl}_t^\ve-\widehat{\qtl}_t|\right)+V\Prob(\tau_1\geq
T(\ve))}
\\
\play{\leq U \max\li_{0\leq t \leq T(\ve)}
\left(\E_{\qtl}|\widehat{\qtl}_t^\ve-\widehat{\qtl}_t|^2\right)^{1/2}
+V\Prob(\tau_1\geq T(\ve))} \ .
\end{array}$$

By the integral form of the stochastic differential equations of the
processes $\widehat{\qtl}_t^\ve$ and $\widehat{\qtl}_t$ we have

$$\begin{array}{l}
|\widehat{\qtl}_t^\ve-\widehat{\qtl}_t|^2
\\
\play{\leq C\left(\left|\int_0^t
\left[\left(\dfrac{B(\widehat{\qtl}_s^\ve)}{\widehat{\Lb}(\widehat{\qtl}_s^\ve)+\ve}
-\dfrac{\widehat{\Lb}'(\widehat{\qtl}_s^\ve)}{2(\widehat{\Lb}(\widehat{\qtl}_s^\ve)+\ve)^3}\right)
-\left(\dfrac{B(\widehat{\qtl}_s^\ve)}{\widehat{\Lb}(\widehat{\qtl}_s^\ve)}
-\dfrac{\widehat{\Lb}'(\widehat{\qtl}_s^\ve)}{2(\widehat{\Lb}(\widehat{\qtl}_s^\ve))^3}\right)\right]ds\right|^2\right.
+}
\\
\play{\left|\int_0^t
\left[\left(\dfrac{B(\widehat{\qtl}_s^\ve)}{\widehat{\Lb}(\widehat{\qtl}_s^\ve)}
-\dfrac{\widehat{\Lb}'(\widehat{\qtl}_s^\ve)}{2(\widehat{\Lb}(\widehat{\qtl}_s^\ve))^3}\right)
-\left(\dfrac{B(\widehat{\qtl}_s)}{\widehat{\Lb}(\widehat{\qtl}_s)}
-\dfrac{\widehat{\Lb}'(\widehat{\qtl}_s)}{2(\widehat{\Lb}(\widehat{\qtl}_s))^3}\right)\right]ds\right|^2
+}
\\
\play{\left.\left|\int_0^t
\left[\dfrac{1}{\widehat{\Lb}(\widehat{\qtl}_s^\ve)+\ve}-\dfrac{1}{\widehat{\Lb}(\widehat{\qtl}_s^\ve)}\right]dW_s\right|^2
+\left|\int_0^t
\left[\dfrac{1}{\widehat{\Lb}(\widehat{\qtl}_s^\ve)}-\dfrac{1}{\widehat{\Lb}(\widehat{\qtl}_s)}\right]dW_s\right|^2\right)}
\ .
\end{array}$$

Let $\al(\lb)$ be the Lipschitz constant of $\dfrac{1}{x}$
($x>\lb$), $\bt(\lb)$ that of $\dfrac{1}{2x^3}$ ($x>\lb$),
$\gm(\dt')$ that of
$\dfrac{B(\widehat{q})}{\widehat{\Lb}(q)}-\dfrac{\widehat{\Lb}'(q)}{2\widehat{\Lb}(q)^3}$
$(q\geq \dt')$, $\mu(\dt')$ that of $\dfrac{1}{\widehat{\Lb}(q)}$
($q\geq \dt'$). Let $m(\dt')\equiv\min\li_{x\in [\dt',b]}\Lb(x)$.

We can estimate

$$\begin{array}{l}
\play{\E_{\qtl}\left|\int_0^t
\left[\left(\dfrac{B(\widehat{\qtl}_s^\ve)}{\widehat{\Lb}(\widehat{\qtl}_s^\ve)+\ve}
-\dfrac{\widehat{\Lb}'(\widehat{\qtl}_s^\ve)}{2(\widehat{\Lb}(\widehat{\qtl}_s^\ve)+\ve)^3}\right)
-\left(\dfrac{B(\widehat{\qtl}_s^\ve)}{\widehat{\Lb}(\widehat{\qtl}_s^\ve)}
-\dfrac{\widehat{\Lb}'(\widehat{\qtl}_s^\ve)}{2(\widehat{\Lb}(\widehat{\qtl}_s^\ve))^3}\right)\right]ds\right|^2}
\\
\play{\leq A_1 (t^2\ve^2[\al^2(m(\dt'))+\bt^2(m(\dt'))]) \ ,}
\end{array}$$

$$\begin{array}{l}
\play{\E_{\qtl}\left|\int_0^t
\left[\left(\dfrac{B(\widehat{\qtl}_s^\ve)}{\widehat{\Lb}(\widehat{\qtl}_s^\ve)}
-\dfrac{\widehat{\Lb}'(\widehat{\qtl}_s^\ve)}{2(\widehat{\Lb}(\widehat{\qtl}_s^\ve))^3}\right)
-\left(\dfrac{B(\widehat{\qtl}_s)}{\widehat{\Lb}(\widehat{\qtl}_s)}
-\dfrac{\widehat{\Lb}'(\widehat{\qtl}_s)}{2(\widehat{\Lb}(\widehat{\qtl}_s))^3}\right)\right]ds\right|^2}
\\
\play{\leq A_2 t \gm^2(\dt')\int_0^t
\E_{\qtl}|\widehat{\qtl}_s^\ve-\widehat{\qtl}_s|^2 ds \ ,}
\end{array}$$

$$\E_{\qtl}\left|\int_0^t
\left[\dfrac{1}{\widehat{\Lb}(\widehat{\qtl}_s^\ve)+\ve}-\dfrac{1}{\widehat{\Lb}(\widehat{\qtl}_s^\ve)}\right]dW_s\right|^2
\leq \int_0^t \ve^2 \al^2(m(\dt'))ds=\ve^2 t \al^2(m(\dt')) \ ,$$

$$\E_{\qtl}\left|\int_0^t
\left[\dfrac{1}{\widehat{\Lb}(\widehat{\qtl}_s^\ve)}-\dfrac{1}{\widehat{\Lb}(\widehat{\qtl}_s)}\right]dW_s\right|^2
\leq \mu^2(\dt')\int_0^t
\E_{\qtl}|\widehat{\qtl}_s^\ve-\widehat{\qtl}_s|^2 ds \ .$$

We have, by using the above estimates, with a possible change of the
constant $C$, that

$$\begin{array}{l}
\play{\E_{\qtl}|\widehat{\qtl}_t^\ve-\widehat{\qtl}_t|^2\leq
C\left(t\ve^2(t(\al^2(m(\dt'))+\bt^2(m(\dt')))+\al^2(m(\dt')))+(t\gm^2(\dt')+\mu^2(\dt'))\int_0^t
\E_{\qtl}|\widehat{\widetilde{q}}_s^\ve-\widehat{\widetilde{q}}_s|^2ds\right)}
\ .
\end{array}$$

By Bellman-Gronwall inequality we have

$$\E_{\qtl}|\widehat{\qtl}_t^\ve-\widehat{\qtl}_t|^2\leq
Ct\ve^2(t(\al^2(m(\dt'))+\bt^2(m(\dt')))
+\al^2(m(\dt')))\exp\left(C(t\gm^2(\dt')+\mu^2(\dt'))t\right) \ .$$

As we can check that $|\al(m(\dt'))|\leq \dfrac{1}{m^2(\dt')}$,
$\bt(m(\dt'))\leq \dfrac{A_3}{m^4(\dt')}$, $\gm(\dt')\leq
\dfrac{A_3}{m^4(\dt')}$ and $|\mu(\dt')|\leq
\dfrac{A_3}{m^2(\dt')}$, this gives, as $\dt'$ is small, that

$$\begin{array}{l}
\max\li_{0\leq t\leq
T(\ve)}\left(\E_{\qtl}|\widehat{\qtl}_t^\ve-\widehat{\qtl}_t|^2\right)^{1/2}\leq
\\
\ \ \ \ \ \ \ \ \ \ \ \leq CT(\ve)\ve(\al^2(m(\dt'))+\bt^2(m(\dt'))
+\dfrac{\al^2(m(\dt'))}{T(\ve)})^{1/2}\exp\left(C(T(\ve)\gm^2(\dt')+\mu^2(\dt'))T(\ve)\right)
\\
\ \ \ \ \ \ \ \ \ \ \ \leq CT(\ve)\dfrac{\ve}{\min\li_{q\in
[1+\dt',1+b]}\lb^4(q)} \exp\left(CT^2(\ve)\dfrac{1}{\min\li_{q\in
[1+\dt',1+b]}\lb^8(q)}\right)
 \ .
\end{array}$$

Noticing that by strong Markov property $\Prob(\tau_1\geq
T(\ve))\leq K\exp(-pT(\ve))$ for some $p>0, K>0$, we see that

$$|\psi_2^\ve(q)|\leq CT(\ve)\dfrac{\ve}{\min\li_{q\in
[1+\dt',1+b]}\lb^4(q)} \exp\left(CT^2(\ve)\dfrac{1}{\min\li_{q\in
[1+\dt',1+b]}\lb^8(q)}\right)+V\exp(-pT(\ve))\ .$$

Let us choose $T(\ve)=\sqrt{\ln\ln\dfrac{1}{\ve}}$. We will then
have

$$|\psi_2^\ve(q)|\leq C\left(\ln\ln\dfrac{1}{\ve}\right)^{1/2}\dfrac{\ve}{\min\li_{q\in
[1+\dt',1+b]}\lb^4(q)} \left(\ln
\dfrac{1}{\ve}\right)^{\frac{C}{\min\li_{q\in
[1+\dt',1+b]}\lb^8(q)}}+V\exp(-p\sqrt{\ln\ln\dfrac{1}{\ve}}))\ .$$

For fixed $\dt'>0$, one can choose $\ve$ small enough such that
$$|\psi_2^\ve(q)|\leq \dfrac{U_0\ve^\kp}{\min\li_{q\in
[1+\dt',1+b]\cup [-1+a,
-1-\dt']}\lb^4(q)}+U_0\exp(-p\sqrt{\ln\ln\dfrac{1}{\ve}})$$ for some
$U_0>0$, $p>0$ and $0<\kp<1$. As we choose first $\ve \da 0$ and
then $\dt'\da 0$, this gives that as $\ve$ is small we have
$|\psi_2^\ve(q)|\leq (\utl(\dt)-\utl(-\dt))^2$. $\square$

\

\textbf{Lemma 2.7.} \textit{We have, as $\ve, \dt, \dt'$ are small,
for $q\in [-1-\dt',1+\dt']$ and $C>0$, that}

$$\left|\Prob_{q}(\pi(q_{\sm_0}^\ve)=\dt)
-\dfrac{\utl(0)-\utl(-\dt)}{\utl(\dt)-\utl(-\dt)}\right|\leq
\dfrac{\utl(\dt')-\utl(0)+C\ve}{\utl(\dt)-\utl(-\dt)} \ , $$

$$\left|\Prob_{q}(\pi(q_{\sm_0}^\ve)=-\dt)
-\dfrac{\utl(\dt)-\utl(0)}{\utl(\dt)-\utl(-\dt)}\right|\leq
\dfrac{\utl(\dt')-\utl(0)+C\ve}{\utl(\dt)-\utl(-\dt)} \ . $$

\

\textbf{Proof.} Let $\qtl=\pi(q)\in [-\dt',\dt']$. We have, for
bounded positive functions $C_1(\dt,\ve)$, $C_2(\dt,\ve)$ and
positive constants $C_1$, $C_2$, $C$, that

$$\begin{array}{l}
\play{\left|\Prob_{q}(\pi(q_{\sm_0}^\ve)=\dt)-
\dfrac{\utl(0)-\utl(-\dt)}{\utl(\dt)-\utl(-\dt)}\right|}
\\
\play{=\left|\dfrac{u^{\ve}(q)-u^\ve(-1-\dt)}{u^\ve(1+\dt)-u^\ve(-1-\dt)}
-\dfrac{\utl(0)-\utl(-\dt)}{\utl(\dt)-\utl(-\dt)}\right|}
\\
\play{=\left|\dfrac{\utl(0)-\utl(-\dt)+\utl(\qtl)-\utl(0)+C_1(\dt,\ve)\ve}{\utl(\dt)-\utl(-\dt)
+C_2(\dt,\ve)\ve}
-\dfrac{\utl(0)-\utl(-\dt)}{\utl(\dt)-\utl(-\dt)}\right|}
\\
\play{\leq\dfrac{(\utl(\qtl)-\utl(0)+C_1\ve)(\utl(\dt)-\utl(-\dt))
+C_2\ve(\utl(0)-\utl(-\dt))}{(\utl(\dt)-\utl(-\dt))^2}}
\\
\play{\leq \dfrac{\utl(\dt')-\utl(0)+C\ve}{\utl(\dt)-\utl(-\dt)} \
.}
\end{array}$$

The estimate of $\Prob_{q}(\pi(q_{\sm_0}^\ve)=-\dt)$ is similar.
$\square$

\

\textbf{Lemma 2.8.} \textit{We have, as $\ve$ are small, for $q\in
[-1-\dt',1+\dt']$, that $|(I)^\ve(q)|\leq C
(\utl(\dt)-\utl(-\dt))^2$. }

\

\textbf{Proof.} We have, using Lemma 2.7, that

$$\begin{array}{l}
|(I)^\ve(q)|
\\
=|\E_qf(\pi(q_{\sm_0}^\ve))-f(\pi(q))|
\\
=|(f(\dt)-f(0))\Prob_q(\pi(q_{\sm_0}^\ve)=\dt)
-(f(0)-f(-\dt))\Prob_q(\pi(q_{\sm_0}^\ve)=-\dt)+(f(0)-f(\pi(q)))|
\\
\play{\leq
\left|(f(\dt)-f(0))\dfrac{\utl(0)-\utl(-\dt)}{\utl(\dt)-\utl(-\dt)}
-(f(0)-f(-\dt))\dfrac{\utl(\dt)-\utl(0)}{\utl(\dt)-\utl(-\dt)}\right|+}
\\
\play{\ \ \ \
C_4\dfrac{\utl(\dt')-\utl(0)+M\ve}{\utl(\dt)-\utl(-\dt)}+C_5(\utl(\dt')-\utl(0))}
\\
\play{=
\left|\dfrac{(\utl(0)-\utl(-\dt))(\utl(\dt)-\utl(0))}{\utl(\dt)-\utl(-\dt)}
\left(\dfrac{f(\dt)-f(0)}{\utl(\dt)-\utl(0)}-\dfrac{f(0)-f(-\dt)}{\utl(0)-\utl(-\dt)}\right)\right|+}
\\
\play{\ \ \ \
C_4\dfrac{\utl(\dt')-\utl(0)+M\ve}{\utl(\dt)-\utl(-\dt)}+C_5(\utl(\dt')-\utl(0))}
\\
\play{\leq
C_3(\utl(\dt)-\utl(-\dt))^2+C_4\dfrac{\utl(\dt')-\utl(0)+M\ve}{\utl(\dt)-\utl(-\dt)}+C_5(\utl(\dt')-\utl(0))}
\ .
\end{array}$$

We have used our gluing condition $D_{\utl}^+f(0)=D_{\utl}^-f(0)$.
Now we choose first $\ve\da 0$ then $\dt' \da 0$, we get, as $\ve$
is small, that $|(I)^\ve(q)|\leq C(\utl(\dt)-\utl(-\dt))^2$.
$\square$

\

\textbf{Lemma 2.9.} \textit{As $\ve, \dt, \dt'$ are small, for $q\in
[-1-\dt',1+\dt']$ we have, }

$$\E_{q}\sm_0\leq
C(\utl(\dt)-\utl(-\dt))(\vtl(\dt)-\vtl(-\dt)) \ , \
\E_q(1-e^{-\lb\sm_0})\leq
C(\utl(\dt)-\utl(-\dt))(\vtl(\dt)-\vtl(-\dt))\ .$$

\

\textbf{Proof.} We apply the well known formula for the expected
exit time (see, for example \cite[Chapter VII, Theorem
3.6]{[Revuz-Yor]}) and we have

$$\E_q\sm_0=\int_{-1-\dt}^{1+\dt}G^\ve(q,r)dv^\ve(r) \ ,$$
where the Green function
$$G^\ve(q,r)=
\left\{
\begin{array}{l}
\dfrac{(u^\ve(q)-u^\ve(-1-\dt))(u^\ve(1+\dt)-u^\ve(r))}{u^\ve(1+\dt)-u^\ve(-1-\dt)}
\text{ for } -1-\dt\leq q\leq r\leq 1+\dt \ ,
\\
\dfrac{(u^\ve(r)-u^\ve(-1-\dt))(u^\ve(1+\dt)-u^\ve(q))}{u^\ve(1+\dt)-u^\ve(-1-\dt)}
\text{ for } -1-\dt\leq r\leq q\leq 1+\dt \ ,
\\
0 \text{ otherwise } \ .
\end{array}
\right.$$

Therefore it is easy to estimate

$$\begin{array}{l}
\E_{q}\sm_0
\\
\play{\leq (u^\ve(1+\dt)-u^\ve(-1-\dt))(v^\ve(1+\dt)-v^\ve(-1-\dt))}
\\
\play{\leq(\utl(\dt)-\utl(-\dt)+C_6\ve)(\vtl(\dt)-\vtl(-\dt)+C_7\ve)}
\\
\leq C(\utl(\dt)-\utl(-\dt))(\vtl(\dt)-\vtl(-\dt))
\end{array}$$ as desired.

This helps us to find

$$\E_{q}(1-e^{-\lb\sm_0})=\lb\E_q\left[\int_0^{\sm_0}e^{-\lb s}ds\right]\leq \lb\E_q\sm_0\leq
C(\utl(\dt)-\utl(-\dt))(\vtl(\dt)-\vtl(-\dt))\ .$$ $\square$

\

\textbf{Lemma 2.10.} \textit{For $q\in G(\dt)$ and $\dt$
sufficiently small, we have}
$$\lim\li_{\dt'\da 0}\lim\li_{\ve\da 0}\E_{q}e^{-\lb \tau_1}\leq 1-C(\utl(\dt))\wedge(-\utl(-\dt)) \
.$$

\

\textbf{Proof.} Without loss of generality let $q\in [1+\dt, 1+b]$.
The expected value $M^\ve(q)=\E_{q}e^{-\lb \tau_1}$ is the solution
of the differential equation $D_{v^\ve}D_{u^\ve}M^\ve(q)=\lb
M^\ve(q)$, $M^\ve(1+\dt')=M^\ve(1+b)=1$.

There exist two solutions $f_1^\lb(q)$, $f_2^\lb(q)$ of the equation
$D_{v}D_{u}f=\lb f$ with $f_1^\lb(1)=f_2^\lb(1+b)=1$ and
$f_1^\lb(1+b)=f_2^\lb(1)=0$. The derivatives $D_{u}f_1^\lb(x)$,
$D_{u}f_2^\lb(x)$ are increasing functions, $-\infty<\lim\li_{q\da
1}D_{u}(f_1^\lb+f_2^\lb)(q)<0$, $0<\lim\li_{q\uparrow
1+b}D_{u}(f_1^\lb+f_2^\lb)(q)<\infty$ (see \cite{[Feller]},
\cite{[Mandl]}).

We shall make use of Lemma 2.6. Since $q\in [1+\dt, 1+b]$ we see
that $\sm_0=0$. Lemma 2.6 tells us that, for $k=1,2$, we have

$$\lim\li_{\ve\da 0}\left|\E_q\left[\int_0^{\tau_1}e^{-\lb t}
[\lb f_k^\lb(q_t^\ve)-D_{v}D_{u}f_k^\lb(q_t^\ve)]dt+e^{-\lb
\tau_1}f_k^\lb(q_{\tau_1}^\ve)\right]-f_k^\lb(q)\right|\leq
(\utl(\dt)-\utl(-\dt))^2 \ .$$

Taking into account the definitions of $f_1^\lb, f_2^\lb$ we see
that the above inequality gives

$$\left|\lim\li_{\ve\da 0}\E_q e^{-\lb
\tau_1}f_k^\lb(q_{\tau_1}^\ve) -f_k^\lb(q)\right|\leq
(\utl(\dt)-\utl(-\dt))^2 \ .$$

Since $f_k^\lb(q_{\tau_1}^\ve)=f_k^\lb(1+\dt')$ when
$q_{\tau_1}^\ve=1+\dt'$ and $f_k^\lb(q_{\tau_1}^\ve)=f_k^\lb(1+b)$
when $q_{\tau_1}^\ve=1+b$, we see that for some $K>0$ we have

$$
\left|\lim\li_{\ve \da 0}\E_{q} e^{-\lb \tau_1}
-\dfrac{(f_2^\lb(1+b)-f_2^\lb(1+\dt'))f_1^\lb(q)+(f_1^\lb(1+\dt')-f_1^\lb(1+b))f_2^\lb(q)}
{f_1^\lb(1+\dt')f_2^\lb(1+b)-f_1^\lb(1+b)f_2^\lb(1+\dt')}
\right|\leq K(\utl(\dt)-\utl(-\dt))^2\ .
$$

(The expression
$$\dfrac{(f_2^\lb(1+b)-f_2^\lb(1+\dt'))f_1^\lb(q)+(f_1^\lb(1+\dt')-f_1^\lb(1+b))f_2^\lb(q)}
{f_1^\lb(1+\dt')f_2^\lb(1+b)-f_1^\lb(1+b)f_2^\lb(1+\dt')}$$ is the
solution of the equation $\lb f(q)=D_vD_uf$ with
$f(1+\dt')=f(1+b)=1$.)

This gives

$$\left|\lim\li_{\dt'\da 0}\lim\li_{\ve\da 0}\E_q(1-e^{-\lb
\tau_1})-[1-(f_1^\lb(q)+f_2^\lb(q))]\right|\leq
K(\utl(\dt)-\utl(-\dt))^2 \ .$$

Taking into account that $-\infty<\lim\li_{q\da
1}D_{u}(f_1^\lb+f_2^\lb)(q)<0$, $0<\lim\li_{q\uparrow
1+b}D_{u}(f_1^\lb+f_2^\lb)(q)<\infty$ we see from the above estimate
that $$\lim\li_{\dt' \da 0}\lim\li_{\ve \da 0}\E_{q}(1-e^{-\lb
\tau_1})\geq C(\utl(\dt))$$ for $q\in [1+\dt, 1+b]$ and $\dt$
sufficiently small. The case of $\utl(-\dt)$ is handled in a similar
way. $\square$

\

\section{A two dimensional model problem}

In this section we discuss a two dimensional model problem. We work
with a Smoluchowski-Kramers approximation in the plane $\R^2$. Let
us suppose that the friction coefficient $\lb(\bullet)$ depends on
the $y$ variable only: $\lb(x,y)=\lb(y)$. Suppose for $y\in [-1,1]$
we have $\lb(y)=0$. For $y\not \in [-1,1]$ we have $\lb(y)>0$. For
simplicity of presentation we also assume that the drift is zero:
$\boldb(\bullet)=\mathbf{0}$. All the other assumptions about
$\lb(\bullet)$ are the same as was made in Section 1.

In addition, we assume that for $\ve>0$,
$$\play{\int_{-\ve-1}^{-1}\dfrac{1}{\lb(y)}dy=\int_{1}^{1+\ve}\dfrac{1}{\lb(y)}dy=\infty} \
.$$ (In the case that both integrals converge the proof of Lemma 3.1
repeat that in the case of both integrals divergent but we do not
know anything about the case of one integral convergent and the
other divergent.)

As we already introduced in equation (1.8) of Section 1, we are
actually considering the stochastic differential equation for the
position of the particle $\boldq_t^\ve\in \R^2$ as follows:

$$\dot{\boldq}^\ve_t=-\dfrac{\grad \lb (\boldq_t^\ve)}
{2(\lb
(\boldq_t^\ve)+\ve)^3}+\dfrac{1}{\lb(\boldq_t^\ve)+\ve}\dot{\boldW}_t
\ , \ \boldq_0^\ve=\boldq_0\in \R^2 \ , \  \ve>0 \ . \eqno(3.1)$$

By taking into account our assumption on the friction coefficient
$\lb$ we can write the above equation in coordinate form. Let
$\boldq_t^\ve=(x_t^\ve,y_t^\ve)$. Let $\boldW_t=(W_t^1,W_t^2)$. We
have

$$\left\{\begin{array}{l}
\dot{x}_t^\ve=\dfrac{1}{\lb(y_t^\ve)+\ve}\dot{W}_t^1 \ , \
x_0^\ve=x_0\in \R \ ,
\\
\dot{y}_t^\ve=-\dfrac{\lb'(y_t^\ve)}{2(\lb(y_t^\ve)+\ve)^3}+\dfrac{1}{\lb(y_t^\ve)+\ve}\dot{W}_t^2
\ , \ y_0^\ve=y_0\in\R \ . \end{array}\right. \eqno(3.2)$$

Let $a<0<b$ be given. Throughout this section we will assume that
our process $\boldq_t^\ve$ is stopped once it exits from the domain
$\{(x,y)\in \R^2: a-1\leq y \leq b+1\}$. We therefore suppose that
$y_0\in [a-1,b+1]$.

Note that, similarly as in Section 2, the process $y_t^\ve$ is a
strong Markov process subject to a generalized second order
differential operator in the form $D_{v^\ve(y)}D_{u^\ve(y)}$ where

$$u^\ve(y)=\int_0^y (\lb(s)+\ve)ds \ , \ v^\ve(y)=2\int_0^y (\lb(s)+\ve)ds \ . \eqno(3.3)$$

Let

$$u(y)=\int_0^y \lb(s)ds \ , \ v(y)=2\int_0^y \lb(s)ds \ . \eqno(3.4)$$

We have the obvious relation $u^\ve(y)=u(y)+\ve y$ and
$v^\ve(y)=v(y)+2\ve y$.

\

Let us identify points in the $x$ direction $x\sim x+2\pi$.
Therefore we get a process on the cylinder $S^1\times [a-1,b+1]$,
stopped once it hits the boundary $\{y=a-1 \text{ or } b+1\}$. Let
$$\left\{\begin{array}{l}
\tht_t^\ve=x_t^\ve \mod 2\pi \ ,
\\
y_t^\ve=y_t^\ve \ .
\end{array}\right.$$

In the rest of this section we refer to the process $\boldq_t^\ve$
as the one on a cylinder: $\boldq_t^\ve=(\tht_t^\ve,y_t^\ve)$ is on
the cylinder $S^1\times[a-1,b+1]$. When we speak about the process
$\boldq_t^\ve$ on the domain $\{(x,y)\in \R^2: a-1\leq y \leq
b+1\}\subset \R^2$ we will instead refer to the coordinate
representation $(x_t^\ve, y_t^\ve)$.

Let $\fC$ be the product $S^1\times [a,b]$ with all points
$S^1\times \{0\}$ identified, forming the point $\fo$. A generic
point on $\fC$ will be denoted $\boldqtl=(\tht,\ytl)$ where $\tht\in
S^1$ and $\ytl\in [a,b]$. All points $(\tht, 0)$ correspond to
$\fo$.

Let us consider the following projection map $\boldpi: S^1\times
[a-1,b+1]\ra \fC$. We let

$$\boldpi(\tht,y)=\left\{
\begin{array}{l}
(\tht, y-1) \ , \ \text{ for } 1 < y\leq b+1 \ ;
\\
(\tht, y+1) \ , \ \text{ for } a-1 \leq y < -1 \ ;
\\
\fo \ , \ \text{ for } -1\leq y\leq 1 \ .
\end{array}\right. \eqno(3.5)$$

Let $\boldpi(\boldq_t^\ve)=\boldqtl_t^\ve=(\tht_t^\ve, \ytl_t^\ve)$.
We see that $\ytl_t^\ve=\pi(y_t^\ve)$ where $\pi$ is the projection
map introduced in Section 2.

Let, as in Section 2, $\utl(\ytl)=u(\ytl-1)$ for $\ytl<0$ and
$\utl(\ytl)=u(\ytl+1)$ for $\ytl>0$ and $\utl(0)=u(1)=u(-1)$;
$\vtl(\ytl)=v(\ytl-1)$ for $\ytl<0$ and $\vtl(\ytl)=v(\ytl+1)$ for
$\ytl>0$ and $\vtl(0)=\vtl(1)=\vtl(-1)$. The functions $\utl(\ytl)$
and $\vtl(\ytl)$ are continuous strictly increasing functions on
$[a,b]$. Let $\lbtl(\ytl)=\lb(\ytl-1)$ for $\ytl<0$ and
$\lbtl(\ytl)=\lb(\ytl+1)$ for $\ytl>0$ and $\lbtl(0)=0$.

Let $A$ be the operator given, for $\ytl\neq 0$, by the formula

$$A f(\tht, \ytl)=D_{\utl(\ytl)}D_{\vtl(\ytl)}f+
\dfrac{1}{\lbtl^2(\ytl)}\dfrac{\pt^2}{\pt \tht^2}f \ . \eqno(3.6)$$

Let $D(A)$ be the subset of the space $\contfunc(\fC)$ consisting of
functions $f(\boldqtl)$ for which $A f(\tht,\ytl)$ is defined and
continuous for $\ytl \neq 0$, the derivatives in it being
continuous; such that finite limits
$$\lim\li_{\tht'\ra \tht, \ytl\ra 0-}D_{\utl(\ytl)}f(\tht', \ytl) \ ,
\ \lim\li_{\tht'\ra \tht, \ytl\ra 0+}D_{\utl(\ytl)}f(\tht', \ytl) \
, \eqno(3.7)$$ exist; $$\lim\li_{\tht' \ra \tht, \ytl \ra
0}Af(\tht', \ytl) \eqno(3.8)$$ exists and does not depend on $\tht$;
$$\lim\li_{\tht'\ra \tht, \ytl \ra a}A f(\tht',
\ytl)=\lim\li_{\tht'\ra \tht, \ytl \ra b}A f(\tht', \ytl)=0 \ ;
\eqno(3.9)$$ and $$\int_0^{2\pi}\lim\li_{\tht' \ra \tht, \ytl \ra
0-}D_{\utl(\ytl)}f(\tht', \ytl)d\tht=\int_0^{2\pi}\lim\li_{\tht' \ra
\tht, \ytl \ra 0+}D_{\utl(\ytl)}f(\tht', \ytl)d\tht \ .
\eqno(3.10)$$

It is worth mentioning here that the above condition (3.10) in the
definition of $D(A)$ can be replaced by the condition that
$\lim\li_{\tht' \ra \tht, \ytl \ra 0-}D_{\utl(\ytl)}f(\tht', \ytl)$
and $\lim\li_{\tht' \ra \tht, \ytl \ra 0+}D_{\utl(\ytl)}f(\tht',
\ytl)$ not depending on $\tht$ and coinciding. In this case the
proof of Lemma 3.1 remains the same.

Let us define, for $f\in D(A)$, $Af(\tht, a)$ and $Af(\tht, b)$ as
the limits (3.9) and $Af(\fo)$ as the limit (3.8). The operator $A$
defined on $D(A)$ is a linear operator $D(A) \mapsto
\contfunc(\fC)$.

\

\textbf{Lemma 3.1.} \textit{The closure $\overline{A|_{D(A)}}$ of
the operator $A|_{D(A)}$ exists and is the infinitesimal operator of
a Markov semigroup on $\contfunc(\fC)$.}

(The corresponding Markov process $\boldqtl_t$ stops after reaching
the boundary of $\fC$ ($\ytl=a \text{ or } b$).)

\

\textbf{Proof. } We use the Hille-Yosida theorem and we check the
following:

$\bullet$ The domain $D(A)$ is dense in $\contfunc(\fC)$.

This is because we can approximate every function $g$ in
$\contfunc(\fC)$ by a function $f$ which is smooth, close to $g$
outside a neighborhood of $\fo$ and is equal to $g(\fo)$ in the
neighborhood of $\fo$. This function $f$ satisfies our restrictions
on $D(A)$ and can approximate the function $g$ with respect to the
norm of $\contfunc(\fC)$ as we choose the neighborhood of $\fo$
small enough.

$\bullet$ The operator $A|_{D(A)}$ satisfies the maximum principle:
for $f\in D(A)$, if this function reaches its maximum value at a
point $\boldqtl\in \fC$ we have $A f (\boldqtl)\leq 0$.

Indeed, for $\boldqtl=(\tht,a)$ or $(\tht,b)$, we have
$Af(\boldqtl)=0$. If $\boldqtl=(\tht,\ytl)$, $\ytl \neq 0$ the first
partial derivatives at $\boldqtl$ are equal to $0$ and
$\dfrac{\pt^2}{\pt \tht^2}f(\tht, \ytl)\leq 0$,
$D_{\vtl(\ytl)}D_{\utl(\ytl)}\leq 0$. Finally, if $\boldqtl=\fo$ we
have the left-hand derivative $D^-_{\utl(\ytl)}f(\tht,0)\geq 0$, the
right-hand derivative $D^+_{\utl(\ytl)}f(\tht,0)\leq 0$ and by
(3.10) both these derivatives are equal to $0$. It follows then that
the limit as $\ytl \ra 0$ of the second $\ytl$-derivative is
non-positive for all $\tht \in S^1$. Since the integral over $S^1$
of the second $\tht$ derivative is equal to $0$ for all $\ytl \neq
0$, taking into account that $A f(\fo)$ is equal to the limit (3.8),
we have that $Af(\fo)\leq 0$.

It follows from the maximum principle that for $\lb>0$ the operator
$\lb I- A|_{D(A)}$ does not send to zero any function that is not
equal 0, and this linear operator has an inverse (that is not
defined on the whole $\contfunc(\fC)$), with $\|(\lb I -
A|_{D(A)})^{-1}\|\leq \lb^{-1}$. Every bounded linear operator does
have a closure (which is just its extension by continuity), and with
it the operators $\lb I - A|_{D(A)}$ and $A|_{D(A)}$ also have
closures.

$\bullet$ Finally, to check that we can apply Hille-Yosida theorem
to the closure $\overline{A|_{D(A)}}$ we have only to check that the
bounded operator $(\lb I - A|_{D(A)})^{-1}$ is defined on a dense
set. That is, for a dense subset of $F\in \contfunc(\fC)$ there
exists a solution $f\in D(A)$ of the equation $$\lb f -A f =F \ .
\eqno(3.11)$$

Let us take $F(\tht,\ytl)=e^{in\tht}G(\ytl)$, defining $F(\fo)$ as
its limit as $\ytl \ra 0$. Of course for $n \neq 0$ we have to have
$\lim\li_{\ytl \ra 0}G(\ytl)$ (which limit we'll take as the value
$G(0)$) equal to $0$.

We shall look for the solution $f\in D(A)$ of the equation (3.11) in
the form $f(\tht,\ytl)=e^{in\tht}g(\ytl)$ (again, for $n \neq 0$ it
should be $g(0)=\lim\li_{\ytl \ra 0}g(\ytl)=0$).

The differential equation for $g(\ytl)$ following from (3.11) is the
ordinary differential equation
$$(\lb + \dfrac{n^2}{\lbtl^2(\ytl)})g(\ytl)-
D_{\vtl(\ytl)}D_{\utl(\ytl)}g(\ytl)=G(\ytl)\ , \eqno(3.12)$$ and it
should be solved with the boundary conditions
$\dfrac{n^2}{\lbtl^2(a)}g(a)-D_{\vtl(\ytl)}D_{\utl(\ytl)}g(a)=
\dfrac{n^2}{\lbtl^2(b)}g(b)-D_{\vtl(\ytl)}D_{\utl(\ytl)}g(b)=0$,
$D_{\utl(\ytl)}^- g(0)= D_{\utl(\ytl)}^+ g(0)$ and for $n \neq 0$,
$g(0)=0$. From the boundary conditions we get at once
$g(a)=\lb^{-1}G(a)$ and $g(b)=\lb^{-1}G(b)$.

For $n=0$ the equation (3.12) with the boundary conditions
$D_{\utl(\ytl)}D_{\vtl(\ytl)}g(a)=D_{\utl(\ytl)}D_{\vtl(\ytl)}g(b)=0$
and the gluing condition $D_{\utl(\ytl)}^- g(0)=D_{\utl(\ytl)}^+
g(0)$ is just the ordinary differential equation for a
one-dimensional diffusion process that has been considered
infinitely many times, and it has a solution for every $G\in
\contfunc[a,b]$. Let us go to the case $n \neq 0$. We are going to
consider the intervals $[a,0)$ and $(0,b]$ separately; what follows
is about the interval $(0,b]$.

Similarly to how it is done in, e.g.\cite{[Feller]}, we can prove
that there exist two non-negative solutions $\xi_1(\ytl)$ and
$\xi_2(\ytl)$ of the equation
$$(\lb+\dfrac{n^2}{\lbtl^2(\ytl)})\xi_i(\ytl)-
D_{\vtl(\ytl)}D_{\utl(\ytl)}\xi_i(\ytl)=0 \ , \ 0<\ytl\leq b \ ,
\eqno(3.13)$$ the first one increasing and the second one
decreasing, $\xi_1(0)=\xi_2(b)=0$, $\xi_1(b)<\infty$,
$\xi_2(0+)=\infty$. The derivatives $D_{\utl(\ytl)}\xi_i(\ytl)$ are
increasing, $D_{\utl(\ytl)}\xi_1(0)=0$, $D_{\utl(\ytl)}\xi_2(b)<0$.

It is easily checked that the Wronskian
$$W(\ytl)=\det \begin{pmatrix}D_{\utl(\ytl)}\xi_1(\ytl) & D_{\utl(\ytl)}\xi_2(\ytl) \\
\xi_1(\ytl) & \xi_2(\ytl)
\end{pmatrix}$$
(both summands $D_{\utl(\ytl)}\xi_1(\ytl)\cdot \xi_2(\ytl)$ and
$-D_{\utl(\ytl)}\xi_2(\ytl)\cdot \xi_1(\ytl)$ are positive) does not
depend on $\ytl$: $W(\ytl)\equiv W>0$.

Now we define, for $\ytl\in [0,b]$,
$$\widetilde{g}(\ytl)=\dfrac{1}{W}\left[\xi_2(\ytl)\int_0^{\ytl}\xi_1(z)\cdot G(z) d\vtl(z)
+\xi_1(\ytl)\int_{\ytl}^b \xi_2(z)\cdot G(z) d\vtl(z)\right] \ .
\eqno(3.14)$$ It is easily checked that $\lb \widetilde{g}(\ytl)-A
\widetilde{g}(\ytl)=G(\ytl)$ for $0<\ytl\leq b$.

Of course

$$|\widetilde{g}(\ytl)|\leq \dfrac{\|G\|}{W}
\left[\xi_2(\ytl)\int_0^{\ytl}\xi_1(z)d\vtl(z)+\xi_1(\ytl)\int_{\ytl}^b
\xi_2(z)d\vtl(z)\right] \ . \eqno(3.15)$$ Let us check that this
goes to $0$ as $y \ra 0+$.

We have:
$$\xi_i(z)=\dfrac{D_{\vtl(\ytl)}D_{\utl(\ytl)}\xi_i(z)}{\lb+n^2/\lbtl^2(z)}$$
so the first summand in the brackets in (3.15) is less or equal

$$\xi_2(\ytl)\cdot \dfrac{D_{\utl(\ytl)}\xi_1(\ytl)-D_{\utl(\ytl)}\xi_1(0)}
{\min\li_{0\leq z \leq
\ytl}[\lb+n^2/\lbtl^2(z)]}=\dfrac{\xi_2(\ytl)\cdot
D_{\utl(\ytl)}\xi_1(\ytl)}{\min\li_{0\leq z \leq
\ytl}[\lb+n^2/\lbtl^2(z)]}<\dfrac{W}{\min\li_{0\leq z \leq
\ytl}[\lb+n^2/\lbtl^2(z)]} \ ,  $$ and it goes to zero as $\ytl \ra
0+$.

The second summand in (3.15) is less or equal

$$\xi_1(\ytl)\cdot \dfrac{D_{\utl(\ytl)}\xi_2(c)-D_{\utl(\ytl)}\xi_2(\ytl)}
{\min\li_{\ytl\leq z \leq c}[\lb+n^2/\lbtl^2(z)]}+\xi_1(\ytl)\cdot
\dfrac{D_{\utl(\ytl)}\xi_2(b)-D_{\utl(\ytl)}\xi_2(c)}{\min\li_{c\leq
z \leq b}[\lb+n^2/\lbtl^2(z)]} \ , \eqno(3.16)$$ where $\ytl<c<b$.
The first term in (3.16) is less or equal
$$\dfrac{-\xi_1(y)\cdot D_{\utl(\ytl)}\xi_2(\ytl)}{\min\li_{\ytl\leq z \leq c}[\lb+n^2/\lbtl^2(z)]}
\leq \dfrac{W}{\min\li_{\ytl\leq z \leq c}[\lb+n^2/\lbtl^2(z)]} \ ,
$$ and it can be made arbitrarily small by choosing a
positive $c$ close enough to $0$. The second term in (3.16), for a
fixed $c>0$, converges to $0$ as $\ytl \ra 0+$. So we get that
$\lim\li_{\ytl\ra 0+}\widetilde{g}(\ytl)=0$.

Now we are going to find $D_{\utl(\ytl)}\widetilde{g}(0+)$. We have:

$$D_{\utl(\ytl)}\widetilde{g}(\ytl)=\dfrac{1}{W}
\left[D_{\utl(\ytl)}\xi_1(\ytl)\int_{\ytl}^b \xi_2(z)\cdot
G(z)d\vtl(z) +D_{\utl(\ytl)}\xi_2(\ytl)\int_0^{\ytl} \xi_1(z)\cdot
G(z)d\vtl(z)\right] \ . \eqno(3.17)$$

The first integral here is equal to $\play{\int_{\ytl}^c+\int_c^b}$,
and it is not greater than
$$\|G\|\cdot[\xi_2(\ytl)\cdot \vtl(c)+\xi_2(c)\cdot \vtl(b)] \ ,$$
and the first summand is not greater than $$\|G\|/W\cdot [W \cdot
\vtl(c)+\xi_2(c)\cdot \vtl(b)\cdot D_{\utl(\ytl)}\xi_1(\ytl)] \ .$$

By choosing $c\in (0,b)$ close enough to $0$ we make $\vtl(c)$
arbitrarily small; and we know $D_{\utl(\ytl)}\xi_1(\ytl)\ra 0$ as
$\ytl \ra 0+$. So the first summand in (3.17) goes to $0$ as $\ytl
\ra 0+$.

The second summand in (3.17) does not exceed in absolute value

$$\|G\|\cdot \xi_1(\ytl)\cdot |D_{\utl(\ytl)}\xi_2(\ytl)|\cdot \vtl(\ytl)
\leq \|G\|\cdot W \cdot \vtl(\ytl)\ra 0  \ \ (\ytl \ra 0+) \ .
$$

Now we are looking for the solution $g(\ytl)$ of the equation (3.12)
with the boundary conditions under this formula in the form
$g(\ytl)=\widetilde{g}(\ytl)+C\cdot \xi_1(\ytl)$. For the
undetermined coefficient $C$ we get one linear equation, and it does
have a solution since $\xi_1(b)\neq 0$.

The same way we get, for $n \neq 0$, a solution $g(\ytl)$ for
$\ytl<0$ with $g(0-)=D_{\utl(\ytl)}g(0-)=0$, $g(a)=\mu^{-1}G(a)$.

So we get a solution $f\in D(A)$ of the equation (3.11) for every
function $F(\tht,\ytl)=\sum\li_{n=-N}^N e^{i n \tht}\cdot
G_n(\ytl)$, $G_n(\ytl)\in \contfunc[a,b]$, such that $G_n(0)=0$ for
$n \neq 0$ (we take $f(\fo)=G_0(0)$). The set of such functions is
dense in $\contfunc(\fC)$ so that the closure operator
$\overline{(\lb I -A|_{D(A)})^{-1}}$ is defined on the whole
$\contfunc(\fC)$ which finishes the proof. $\square$

\

Let $\boldqtl_t$ be the Markov process corresponding to
$\overline{A|_{D(A)}}$, whose existence was proved in Lemma 3.1. We
prove the following

\

\textbf{Theorem 3.1.}  \textit{As $\ve \da 0$, for fixed $T>0$, the
process $\boldqtl_t^\ve=\boldpi(\boldq_t^\ve)$ converges weakly in
the space $\mathbf{C}_{[0,T]}(\fC)$ to the process $\boldqtl_t$.}

\

The proof is again based on an application of Lemma 2.2.

\

\textbf{Proof of Theorem 3.1.} Making use of Lemma 2.2, we take the
metric space $M=S^1\times [a-1,b+1]$ with standard metric. The
mapping $Y=\boldpi$. The space $Y(M)=\fC$ is endowed with the metric
$d$, defined as follows. For any two points $(\tht_1,\ytl_1)$ and
$(\tht_2,\ytl_2)$ on $\fC$ with $\ytl_1, \ytl_2$ having the same
sign we let $d((\tht_1,\ytl_1), (\tht_2,\ytl_2))$ be the Euclidean
distance between points $(|\ytl_1|\cos \tht_1, |\ytl_1|\sin \tht_1)$
and $(|\ytl_2|\cos \tht_2, |\ytl_2|\sin \tht_2)$ in $\R^2$; if
$\ytl_1$ and $\ytl_2$ have different sign we take
$d((\tht_1,\ytl_1),
(\tht_2,\ytl_2))=d((\tht_1,\ytl_1),\fo)+d(\fo,(\tht_2,\ytl_2))$.
With respect to this metric the space $\fC$ is a complete separable
metric space. We take the process $(X_t^\ve, \Prob_x^\ve)$ as
$\boldq_t^\ve$ and the process $(y_t, \Prob_y)$ is taken as
$\boldqtl_t$.

For the uniqueness of solution of martingale problem we set the
space $\Psi$ be the space of all continuous functions on $\fC$ which
has the form $F(\tht,\ytl)=\sum\li_{n=-N}^N e^{in\tht}\cdot
G_n(\ytl)$, $G_n\in \contfunc[a,b]$ is continuously differentiable
inside $[a,0)$ and $(0,b]$, also $G_n(0)=0$ for $n \neq 0$. We take
$f(\fo)=G_0(0)$. It is proved in the proof of Lemma 3.1 that the
equation $\lb f- A f= F$ always has a solution $f\in D\subset D(A)$
for all $F\in \Psi$ and $\lb>0$. The space $D$ contains those
functions $f\in \contfunc(\fC)$ that are bounded and are three times
continuously differentiable inside $\fC^+\equiv \{(\tht, \ytl)\in
\fC: a<\ytl<0\}$ and $\fC^-\equiv \{(\tht, \ytl)\in \fC:
0<\ytl<b\}$.

We will state pre-compactness of family of distributions of
processes $\boldqtl_t^\ve$ in Lemma 3.2. What remains to do is to
check that for every compact $K \subset \fC$ and for every $f\in D$
and every $\lb>0$ we have

$$\E_{\boldq_0}\left[\int_0^\infty e^{-\lb t}[\lb f(\boldpi(\boldq_t^\ve))
-A f(\boldpi(\boldq_t^\ve))]dt- f(\boldpi(\boldq_0))\right]\ra 0$$
as $\ve \da 0$ uniformly in $\boldq_0\in \boldpi^{-1}(K)$. The proof
of this is essentially the same as the proof we did in Lemma 2.5,
based on the following auxiliary Lemmas 3.9 (for the proof of
convergence for processes near $\fo$) and 3.10 (for the proof of
convergence for processes away from $\fo$) and the auxiliary Lemmas
2.9 and 2.10 (for the estimates on the exit times, notice that the
stopping times $\sm_n$ and $\tau_n$ we will work with in this
section are essentially the same stopping times that we worked with
in Section 2 since we are discussing about a model problem). We omit
the details in the proof. $\square$

\

Let $\kp$ be a real number with small absolute value. Let
$G(\kp)=\{(\tht,y)\in S^1\times [a-1,b+1]: a-1\leq y\leq -1-\kp
\text{ or } 1+\kp \leq y \leq b+1\}$. Let $C^+(\kp)=\{(\tht,y)\in
S^1\times [a-1,b+1]: y=1+\kp\}$ and $C^-(\kp)=\{(\tht,y)\in
S^1\times [a-1,b+1]: y=-1-\kp\}$. Let $C(\kp)=C^+(\kp)\cup
C^-(\kp)$. Let $\dt>\dt'>0$ be small. We shall introduce a sequence
of stopping times $\tau_0\leq \sm_0<\tau_1<\sm_1<\tau_2<\sm_2<...$
by

$$\tau_0=0 \ , \ \sm_n=\min\{t\geq \tau_n, \boldq_t^\ve\in G(\dt)\} \ ,
\ \tau_n=\min\{t \geq \sm_{n-1}, \boldq_t^\ve\in C(\dt')\} \ .
$$

This is well-defined up to some $\sm_k$ ($k \geq 0$) such that
$$\Prob_{y_{\sm_k}^\ve}(y_{t+\sm_k}^\ve \text{ hits } a-1 \text{ or } b+1
\text{ before it hits } -1-\dt' \text{ or } 1+\dt')=1 \ .$$

We will then define $\tau_{k+1}=\min\{t>\sm_k: y_t^\ve=a-1 \text{ or
} b+1\}$. And we define
$\tau_{k+1}<\sm_{k+1}=\tau_{k+1}+1<\tau_{k+2}=\tau_{k+1}+2<\sm_{k+2}=\tau_{k+1}+3<...$
and so on.

We have $\lim\li_{n \ra \infty}\tau_n=\lim\li_{n \ra
\infty}\sm_n=\infty$. And we have obvious relations
$\boldq_{\tau_n}^\ve\in C(\dt')$, $\boldq_{\sm_n}^\ve\in C(\dt)$ for
$1\leq n \leq k$ (as long as $k \geq 1$, if $k=0$ the process may
start from $G(\dt)$ and goes directly to $S^1\times\{a-1\}$ or
$S^1\times\{b+1\}$ without touching $C(\dt')$ and is stopped there,
or it may start from $S^1\times(-1-\dt,1+\dt)$, reaches $C(\dt)$
first and then goes directly to $S^1\times\{a-1\}$ or
$S^1\times\{b+1\}$ without touching $C(\dt')$ and is stopped there).
Also, for $n\geq k+1$ we have
$\boldq_{\tau_n}^\ve=\boldq_{\sm_n}^\ve\in S^1\times\{a-1\} \text{
or } S^1\times\{b+1\}$. If $\boldq_0^\ve=\boldq_0\in G(\dt)$, then
we have $\sm_0=0$ and $\tau_1$ is the first time at which the
process $\boldq_t^\ve$ reaches $C(\dt')$ or $S^1\times\{a-1\}$ or
$S^1\times\{b+1\}$.

Note that these stopping times are the same as those defined in
Section 2 since our process $y_t^\ve$ is essentially the process
$q_t^\ve$ in Section 2.

\

The pre-compactness of the family $\{\boldqtl_t^\ve\}_{\ve>0}$ in
$\contfunc_{[0,T]}(\fC)$ for $0<T<\infty$ is proved in the same way
as in the one-dimensional case. We shall make use of the technical
Lemma 2.3 with $\qtl_\bullet^{\ve,\dt}$ and $\qtl_\bullet^\ve$
replaced by $\boldqtl_\bullet^{\ve,\dt}$ and $\boldqtl_\bullet^\ve$
and the space $\contfunc_{[0,T]}(\fC)$ instead of
$\contfunc_{[0,T]}([a,b])$. We omit the proof of the next lemma.

\

\textbf{Lemma 3.2.} \textit{The family of distributions of
$\{\boldqtl_t^\ve\}_{\ve>0}$ is pre-compact in
$\contfunc_{[0,T]}(\fC)$.}

\

The next few lemmas establish the estimates on the asymptotic joint
law of the processes $(y_t^\ve,\tht_t^\ve)$ at first exit from a
small neighborhood of the domain within which the friction vanishes.
This is the key part to the proof of Theorem 3.1.

\

Let $\dt''>0$ be small. We consider the process $\boldq_t^\ve$
starting from $\boldq_0^\ve=\boldq_0\in S^1\times [-1-\dt',1+\dt']$.
Let us introduce another sequence of stopping times
$\al_1<\bt_1<\al_2<\bt_2<...<\al_{n(\ve)}$ by
$$\al_1=\min\{0\leq t < \sm_0: \boldq_t^\ve\in C(0)\} \ ,
\ \bt_1=\min\{\al_1< t < \sm_0: \boldq_t^\ve\in C(-\dt'')\} \ ,
$$ and for $k \geq 2$ we define
$$\al_k=\min\{\bt_{k-1}<t<\sm_0: \boldq_t^\ve\in C(0)\} \ , \
\bt_k=\min\{\al_k<t<\sm_0: \boldq_t^\ve\in C(-\dt'')\} \ .
$$

Here we take the convention that the minimum over an empty set is
$\infty$. The number $n(\ve)$ is a non-negative integer-valued
random variable such that $\al_{n(\ve)}<\infty$ and
$\bt_{n(\ve)}=\infty$. If $\al_1=\infty$ we set $n(\ve)=0$.

\

\textbf{Lemma 3.3.} \textit{For $\boldq_0\in G(\dt')$ we have}

$$\Prob_{\boldq_0}(\al_1<\infty)\geq 1-
\max\left(\dfrac{\utl(\dt')+\ve \dt'}{\utl(\dt)+\ve \dt},
\dfrac{-\utl(-\dt')+\ve \dt'}{-\utl(-\dt)+\ve \dt} \right) \ .
\eqno(3.18)$$

\

\textbf{Proof.} If $1\leq y_0^\ve=y_0\leq 1+\dt'$ we have

$$\Prob_{\boldq_0}(\al_1<\infty)=\dfrac{u^\ve(1+\dt)-u^\ve(y)}{u^\ve(1+\dt)-u^\ve(1)}
\geq
\dfrac{u^\ve(1+\dt)-u^\ve(1+\dt')}{u^\ve(1+\dt)-u^\ve(1)}=1-\dfrac{\utl(\dt')+\ve
\dt'}{\utl(\dt)+\ve \dt} \ . $$

If $-1-\dt'\leq y_0^\ve=y_0\leq -1$ we have

$$\Prob_{\boldq_0}(\al_1<\infty)=\dfrac{u^\ve(y)-u^\ve(-1-\dt')}{u^\ve(-1)-u^\ve(-1-\dt)}
\geq
\dfrac{u^\ve(-1-\dt')-u^\ve(-1-\dt)}{u^\ve(-1)-u^\ve(-1-\dt)}=1-\dfrac{-\utl(-\dt')+\ve
\dt'}{-\utl(-\dt)+\ve \dt} \ .$$

If $-1<y^\ve_0=y_0<1$ we have $\Prob_{\boldq_0}(\al_1<\infty)=1$.
$\square$

\

\textbf{Lemma 3.4.} \textit{For $\boldq_0\in G(\dt')$ we have}

$$\Prob_{\boldq_0}(\bt_1<\infty|\al_1<\infty)\geq 1-\max
\left(\dfrac{\ve \dt''}{\utl(\dt)+\ve(\dt+\dt'')}, \dfrac{\ve
\dt''}{-\utl(-\dt)+\ve(\dt+\dt'')}\right)\ . \eqno(3.19)$$

\

\textbf{Proof.} If $y_{\al_1}^\ve=1$ we have

$$\Prob_{\boldq_0}(\bt_1<\infty|\al_1<\infty)=\dfrac{u^\ve(1+\dt)-u^\ve(1)}{u^\ve(1+\dt)-u^\ve(1-\dt'')}
=1-\dfrac{\ve \dt''}{\utl(\dt)+\ve(\dt+\dt'')} \ .$$

If $y_{\al_1}^\ve=-1$ we have

$$\Prob_{\boldq_0}(\bt_1<\infty|\al_1<\infty)=\dfrac{u^\ve(-1)-u^\ve(-1-\dt)}{u^\ve(-1+\dt'')-u^\ve(-1-\dt)}
=1-\dfrac{\ve \dt''}{-\utl(-\dt)+\ve(\dt+\dt'')} \ .$$ $\square$

\

Let $M(\ve)\ra \infty$ as $\ve \da 0$ be an integer. The exact
asymptotics of $M(\ve)$ will be specified later. We prove

\

\textbf{Lemma 3.5.} \textit{For $\boldq_0\in G(\dt')$ we have}

$$\Prob_{\boldq_0}(n(\ve)\geq M(\ve)|\al_1<\infty)\geq \left[1-\max
\left(\dfrac{\ve \dt''}{\utl(\dt)+\ve(\dt+\dt'')}, \dfrac{\ve
\dt''}{-\utl(-\dt)+\ve(\dt+\dt'')}\right)\right]^{M(\ve)-1} \ .
\eqno(3.20)$$

\

\textbf{Proof.} This is because trajectories of $\boldq_t^\ve$
between times $\al_i\leq t<\al_{i+1}$ are independent and by
iteratively using Lemma 3.4 we get the desired result. $\square$

\

\textbf{Lemma 3.6.} \textit{We have} $$\al_{i+1}-\bt_i\geq
\ve^2\left(\dfrac{\dt''}{H_i}\right)^5 \eqno(3.21)$$ \textit{with
$H_i$ being i.i.d. positive random variables with $\E(H_i)^4<\infty$
for $i=1,2,...,n(\ve)-1$.}

\

\textbf{Proof.} This is a result of the H\"{o}lder continuity of the
standard Wiener trajectory $|W_t-W_s|\leq H_i |t-s|^{1/5}$ and the
fact that between times $\bt_i\leq t <\al_{i+1}$ the process
$y_t^\ve$ is a time-changed Wiener process $\dfrac{1}{\ve} W_t$
traveling at least a distance of $\dt''$. $\square$

\

Let us define an auxiliary function

$$\begin{array}{l}
\Om(\ve, \dt, \dt', \dt'', M(\ve))
\\
\equiv 2\left[1-\left[1-\max \left(\dfrac{\ve
\dt''}{\utl(\dt)+\ve(\dt+\dt'')}, \dfrac{\ve
\dt''}{-\utl(-\dt)+\ve(\dt+\dt'')}\right)\right]^{M(\ve)-1}+\right.
\\
\left. \ \ \ \ \ \ \ \ \ \ \ \ 2\max\left(\dfrac{\utl(\dt')+\ve
\dt'}{\utl(\dt)+\ve \dt}, \dfrac{-\utl(-\dt')+\ve
\dt'}{-\utl(-\dt)+\ve \dt} \right) \right] \ .
\end{array}$$

\textbf{Lemma 3.7.} \textit{For $\boldq_0\in G(\dt')$ and for some
$A>0$, $\kp>0$ and $C>0$, there exists $\ve_0>0$ such that for all
$0<\ve<\ve_0$, for any $0\leq \tht_1\leq \tht_2 \leq 2\pi$ we have}

$$\begin{array}{l}
\left|\Prob_{\boldq_0}(\tht_{\sm_0}^\ve\in [\tht_1,\tht_2],
y_{\sm_0}^\ve=1+\dt)-\dfrac{\tht_2-\tht_1}{2\pi}\Prob_{\boldq_0}(y_{\sm_0}^\ve=1+\dt)\right|
\\
\play{\leq C \exp(-A (\dt'')^5 \kp M(\ve)) + 2
\Om(\ve,\dt,\dt',\dt'',M(\ve))}
\end{array}$$
\textit{and}
$$\begin{array}{l}
\left|\Prob_{\boldq_0}(\tht_{\sm_0}^\ve\in [\tht_1,\tht_2],
y_{\sm_0}^\ve=-1-\dt)-\dfrac{\tht_2-\tht_1}{2\pi}\Prob_{\boldq_0}(y_{\sm_0}^\ve=-1-\dt)\right|
\\
\play{\leq C \exp(-A (\dt'')^5 \kp M(\ve)) + 2
\Om(\ve,\dt,\dt',\dt'',M(\ve))} \ .
\end{array}$$

\

\textbf{Proof.} As we have

$$x_t^\ve=\int_0^t\dfrac{1}{\lb(y_t^\ve)+\ve}dW_s^1=
W^1\left(\int_0^t \dfrac{ds}{(\lb(y_s^\ve)+\ve)^2}\right) \ ,$$ we
set $\play{T^\ve(t)=\int_0^t \dfrac{ds}{(\lb(y_s^\ve)+\ve)^2}}$.
Using Lemma 3.6 for $\boldq_0\in G(\dt')$ the random time
$T^\ve(\sm_0)$ can be estimated from below by

$$T^\ve(\sm_0)\geq \int_0^{\sm_0} \dfrac{ds}{(\lb(y_s^\ve)+\ve)^2}\geq
\dfrac{1}{\ve^2}\int_0^{\sm_0}\1_{\{-1\leq y_s^\ve\leq 1\}}ds\geq
\dfrac{1}{\ve^2}\sum\li_{i=1}^{n(\ve)-1}(\al_{i+1}-\bt_i)\geq
(\dt'')^5\sum\li_{i=1}^{n(\ve)-1}\dfrac{1}{(H_i)^5} \ .$$ (If
$n(\ve)=0,1$ the sum is supposed to be $0$.)

And we also notice that the random time $T^\ve(\sm_0)$ only depends
on the behavior of the process $y_t^\ve$ and is therefore
independent of the Wiener process $W_t^1$ in the stochastic
differential equation
$\play{\dot{x}_t^\ve=\dfrac{1}{\lb(y_t^\ve)+\ve}\dot{W}_t^1}$ (see
(3.2)). For the same reason the random variables $y_{\sm_0}^\ve$,
$n(\ve)$ and $\al_1$ are of course also independent of $W_t^1$.

As we have the elementary inequality $\play{\left(\E
\dfrac{1}{(H_i)^5}\right)^{1/5}(\E (H_i)^4)^{1/4}}\geq
\left(\E\dfrac{1}{H_i}\right)(\E H_i)\geq 1$, we have, by Strong Law
of Large Numbers
$$\lim\li_{\ve \da
0}\dfrac{1}{M(\ve)-1}\sum\li_{i=1}^{M(\ve)-1}\dfrac{1}{(H_i)^5}
=\E\left(\dfrac{1}{(H_i)^5}\right)\geq \dfrac{1}{(\E
(H_i)^4)^{5/4}}\geq c >0 \ \text{ a. s. }$$ for some constant $c>0$.
(We can always assume that $H_i$ is uniformly bounded from below by
a positive constant so that
$\left(\E\dfrac{1}{(H_i)^5}\right)<\infty$ and we can apply SLLN.)

Now we see that we can find some $\ve_0>0$ such that for all
$0<\ve<\ve_0$ we will have
$$\Prob_{\boldq_0}(T^\ve(\sm_0)\geq (\dt'')^5\kp M(\ve)|n(\ve)\geq M(\ve), \al_1<\infty)=1 $$
for some constant $\kp>0$.

This gives

$$\begin{array}{l}
\Prob_{\boldq_0}(T^\ve(\sm_0)\geq (\dt'')^5\kp M(\ve),
y_{\sm_0}^\ve=1+\dt|n(\ve)\geq M(\ve), \al_1<\infty)
\\
=\Prob_{\boldq_0}(y_{\sm_0}^\ve=1+\dt|n(\ve)\geq M(\ve),
\al_1<\infty) \ .
\end{array}$$

Recall that we have $\tht_{\sm_0}^\ve= x_{\sm_0}^\ve \text{ mod }
2\pi=W_{T^\ve(\sm_0)}^1 \text{ mod } 2\pi$. Using this, the
independence of $T^\ve(\sm_0)$, $y_{\sm_0}^\ve$, $\al_1$ and
$n(\ve)$ with $W_t^1$, and the above estimates we have, as
$0<\ve<\ve_0$, that

$$\begin{array}{l}
\Prob_{\boldq_0}(\tht_{\sm_0}^\ve\in [\tht_1,\tht_2], y_{\sm_0}^\ve
= 1+ \dt |n(\ve)\geq M(\ve), \al_1<\infty)
\\
\play{= \int_0^\infty \Prob_{\boldq_0}(T^\ve(\sm_0)\in dt,
y_{\sm_0}^\ve=1+\dt|n(\ve)\geq M(\ve),
\al_1<\infty)\Prob_{\boldq_0}(W_t^1 \mod 2 \pi \in [\tht_1,\tht_2])}
\\
\play{= \int_{(\dt'')^5 \lb M(\ve)}^\infty
\Prob_{\boldq_0}(T^\ve(\sm_0)\in dt, y^\ve_{\sm_0}=1+\dt|n(\ve)\geq
M(\ve), \al_1<\infty)\Prob_{\boldq_0}(W_t^1 \mod 2 \pi \in
[\tht_1,\tht_2]) \ .}
\end{array}$$

Since we have the exponential decay $$\left|\Prob(W_t^1\mod 2\pi \in
[\tht_1,\tht_2])-\dfrac{\tht_2-\tht_1}{2\pi}\right|<C\exp(-A t) \
$$ for some $C>0$ and $A>0$, we could estimate
$$\begin{array}{l}
\play{\left|\Prob_{\boldq_0}(\tht_{\sm_0}^\ve\in [\tht_1,\tht_2] \ ,
\ y_{\sm_0}^\ve=1+\dt|n(\ve)\geq M(\ve), \al_1<\infty)-\right.}
\\
\play{\ \ \ \ \ \ \ \ \ \ \ \
\dfrac{\tht_2-\tht_1}{2\pi}\left.\Prob_{\boldq_0}(y^\ve_{\sm_0}=1+\dt|n(\ve)\geq
M(\ve), \al_1<\infty)\right|}
\\
\play{<C\exp(-A (\dt'')^5 \kp M(\ve))}
\\
\end{array}$$ for $0<\ve<\ve_0$.

Notice that we have, by using Lemmas 3.5 and 3.3,

$$\begin{array}{l}
\left|\Prob_{\boldq_0}(\tht_{\sm_0}^\ve \in [\tht_1,\tht_2],
y_{\sm_0}^\ve=1+\dt)-\Prob_{\boldq_0}(\tht_{\sm_0}^\ve\in
[\tht_1,\tht_2] \ , \ y_{\sm_0}^\ve=1+\dt|n(\ve)\geq M(\ve),
\al_1<\infty)\right|
\\
=\left|\Prob_{\boldq_0}(\tht_{\sm_0}^\ve \in [\tht_1,\tht_2],
y_{\sm_0}^\ve=1+\dt|n(\ve)\geq M(\ve), \al_1<\infty)\Prob(n(\ve)\geq
M(\ve), \al_1<\infty)\right.-
\\
\ \ \ \ \ \ \ \left.  \Prob_{\boldq_0}(\tht_{\sm_0}^\ve\in
[\tht_1,\tht_2] \ , \ y_{\sm_0}^\ve=1+\dt|n(\ve)\geq M(\ve),
\al_1<\infty)
\right|+\Prob_{\boldq_0}(n(\ve)<M(\ve))+\Prob_{\boldq_0}(\al_1=\infty)
\\
\leq
2(\Prob_{\boldq_0}(n(\ve)<M(\ve))+\Prob_{\boldq_0}(\al_1=\infty))
\\
\leq
2(\Prob_{\boldq_0}(n(\ve)<M(\ve)|\al_1<\infty)+2\Prob_{\boldq_0}(\al_1=\infty))
\\
\leq 2\left[1-\left[1-\max \left(\dfrac{\ve
\dt''}{\utl(\dt)+\ve(\dt+\dt'')}, \dfrac{\ve
\dt''}{-\utl(-\dt)+\ve(\dt+\dt'')}\right)\right]^{M(\ve)-1}\right.+
\\
\left. \ \ \ \ \ \ \ \ \ \ \ \  2\max\left(\dfrac{\utl(\dt')+\ve
\dt'}{\utl(\dt)+\ve \dt}, \dfrac{-\utl(-\dt')+\ve
\dt'}{-\utl(-\dt)+\ve \dt} \right) \right]
\\
= \Om(\ve, \dt, \dt', \dt'', M) \ .
\end{array}$$

By the same argument we can estimate

$$
\left|\dfrac{\tht_2-\tht_1}{2\pi}\Prob_{\boldq_0}(y_{\sm_0}^\ve=1+\dt)-\dfrac{\tht_2-\tht_1}{2\pi}\Prob_{\boldq_0}(y_{\sm_0}^\ve=1+\dt|n(\ve)\geq
M(\ve), \al_1<\infty)\right|\leq \Om(\ve, \dt, \dt', \dt'', M) \ .
$$

Summing up these estimates we have

$$\begin{array}{l}
\left|\Prob_{\boldq_0}(\tht_{\sm_0}^\ve\in [\tht_1,\tht_2],
y_{\sm_0}^\ve=1+\dt)-\dfrac{\tht_2-\tht_1}{2\pi}\Prob_{\boldq_0}(y_{\sm_0}^\ve=1+\dt)\right|
\\
\leq \left|\Prob_{\boldq_0}(\tht_{\sm_0}^\ve\in [\tht_1,\tht_2],
y_{\sm_0}^\ve=1+\dt)-\Prob_{\boldq_0}(\tht_{\sm_0}^\ve\in
[\tht_1,\tht_2] \ , \ y_{\sm_0}^\ve=1+\dt|n(\ve)\geq M(\ve),
\al_1<\infty)\right|+
\\
\ \ \  \play{\left|\Prob_{\boldq_0}(\tht_{\sm_0}^\ve\in
[\tht_1,\tht_2] \ , \ y_{\sm_0}^\ve=1+\dt|n(\ve)\geq M(\ve),
\al_1<\infty)- \right.}
\\
\play{ \ \ \ \ \ \ \ \ \ \ \ \ \ \ \ \ \ \ \ \
\dfrac{\tht_2-\tht_1}{2\pi}\left.\Prob_{\boldq_0}(y^\ve_{\sm_0}=1+\dt|n(\ve)\geq
M(\ve), \al_1<\infty)\right|} +
\\
\ \ \
\left|\dfrac{\tht_2-\tht_1}{2\pi}\Prob_{\boldq_0}(y_{\sm_0}^\ve=1+\dt)-\dfrac{\tht_2-\tht_1}{2\pi}\Prob_{\boldq_0}(y_{\sm_0}^\ve=1+\dt|n(\ve)\geq
M(\ve), \al_1<\infty)\right|
\\
\leq 2\Om(\ve,\dt,\dt',\dt'',M)+ C \exp(-A(\dt'')^5 \kp M(\ve)) \ ,
\end{array}$$ as desired. The other inequality is established in a similar way. $\square$

\

Combining Lemma 3.7 and Lemma 2.7 we can have

\

\textbf{Lemma 3.8.}  \textit{For $\boldq_0\in G(\dt')$ and for some
$A>0$, $\kp>0$ and $C_1,C_2>0$, there exists $\ve_0>0$ such that for
all $0<\ve<\ve_0$, for any $0\leq \tht_1\leq \tht_2 \leq 2\pi$ we
have}

$$\begin{array}{l}
\left|\Prob_{\boldq_0}(\tht_{\sm_0}^\ve\in [\tht_1,\tht_2],
y_{\sm_0}^\ve=1+\dt)-\dfrac{\tht_2-\tht_1}{2\pi}\dfrac{\utl(0)-\utl(-\dt)}{\utl(\dt)-\utl(-\dt)}\right|
\\
\play{\leq C_1 \exp(-A (\dt'')^5 \kp M(\ve)) + 2
\Om(\ve,\dt,\dt',\dt'',M(\ve)) +
\dfrac{\utl(\dt')-\utl(0)+C_2\ve}{\utl(\dt)-\utl(-\dt)}\equiv
\rho(\ve)} \ ,
\end{array} \eqno(3.22)$$
\textit{and}
$$\begin{array}{l}
\left|\Prob_{\boldq_0}(\tht_{\sm_0}^\ve\in [\tht_1,\tht_2],
y_{\sm_0}^\ve=-1-\dt)-\dfrac{\tht_2-\tht_1}{2\pi}\dfrac{\utl(\dt)-\utl(0)}{\utl(\dt)-\utl(-\dt)}\right|
\\
\play{\leq C_1 \exp(-A (\dt'')^5 \kp M(\ve)) + 2
\Om(\ve,\dt,\dt',\dt'',M(\ve)) +
\dfrac{\utl(\dt')-\utl(0)+C_2\ve}{\utl(\dt)-\utl(-\dt)}\equiv
\rho(\ve)} \ .
\end{array} \eqno(3.23)$$

\

Now let us specify the asymptotic order of $M(\ve)\ra \infty$,
$\dt=\dt(\ve)\ra 0$, $\dt'=\dt'(\ve)\ra 0$ and $\dt''=\dt''(\ve)\ra
0$ as $\ve \da 0$. Since for $0<\kp<1$ we have the elementary
estimate $1-(1-\kp)^n=\kp(1+(1-\kp)+...+(1-\kp)^{n-1})\leq \kp n$ we
can estimate

$$\begin{array}{l}
\Om(\ve, \dt, \dt', \dt'', M(\ve))
\\
\leq 2\left[M(\ve)\cdot\max \left(\dfrac{\ve
\dt''}{\utl(\dt)+\ve(\dt+\dt'')}, \dfrac{\ve
\dt''}{-\utl(-\dt)+\ve(\dt+\dt'')}\right)+\right.
\\
\left. \ \ \ \ \ \ \ \ \ \ \ \  2\max\left(\dfrac{\utl(\dt')+\ve
\dt'}{\utl(\dt)+\ve \dt}, \dfrac{-\utl(-\dt')+\ve
\dt'}{-\utl(-\dt)+\ve \dt} \right) \right] \ .
\end{array}$$

We shall choose $\dt''=\dt''(\ve)<<\dt$ and $M(\ve)$ such that the
requirements of Lemmas 2.6, 2.7 and 2.8 hold. At the same time, we
need
$$(\dt'')^5 M(\ve)\gtrsim \ln\dfrac{1}{(\utl(\dt)-\utl(-\dt))^2} \eqno(3.24)$$ and $$M(\ve)\dfrac{\ve
\dt''}{\utl(\dt)\wedge (-\utl(-\dt))}\lesssim
(\utl(\dt)-\utl(-\dt))^2 \ . \eqno(3.25)$$

To this end we let $M(\ve)=\ln\left(\dfrac{1}{\ve}\right)$ and
$\dt''=\left(\dfrac{(\frac{1}{\utl(\dt)-\utl(-\dt)})\ln(\frac{1}{\utl(\dt)-\utl(-\dt)})^2}{\ln(\frac{1}{\ve})}\right)^{1/5}$.
At the same time we keep our asymptotic order of choice of $\ve$,
$\dt$ and $\dt'$ as in Section 2. This means that we need $$\ve
\left(\ln \left(\dfrac{1}{\ve}\right)\right)^{4/5}
\dfrac{1}{\utl(\dt)-\utl(-\dt)}\ln
\left(\dfrac{1}{\utl(\dt)-\utl(-\dt)}\right)^2\lesssim
(\utl(\dt)-\utl(-\dt))^2 \ .$$ It could be checked that this is
possible to make (3.24) and (3.25) to hold. We formulate this as a
corollary.

\

\textbf{Corollary 3.1.} \textit{Let $\boldq_0\in G(\dt')$. Under the
above specified asymptotic order we have, there exist $\ve_0>0$ such
that for any $0<\ve<\ve_0$ we have}
$$
\left|\Prob_{\boldq_0}(\tht_{\sm_0}^\ve\in [\tht_1,\tht_2],
y_{\sm_0}^\ve=1+\dt)-\dfrac{\tht_2-\tht_1}{2\pi}\dfrac{\utl(0)-\utl(-\dt)}{\utl(\dt)-\utl(-\dt)}\right|
\leq  C\cdot (\utl(\dt)-\utl(-\dt))^2 \ ,
 \eqno(3.26)$$

$$
\left|\Prob_{\boldq_0}(\tht_{\sm_0}^\ve\in [\tht_1,\tht_2],
y_{\sm_0}^\ve=-1-\dt)-\dfrac{\tht_2-\tht_1}{2\pi}\dfrac{\utl(\dt)-\utl(0)}{\utl(\dt)-\utl(-\dt)}\right|
\leq C\cdot (\utl(\dt)-\utl(-\dt))^2 \ .
 \eqno(3.27)$$

\

\textbf{Lemma 3.9.} \textit{For any $\boldq\in G(\dt')$ and for any
$\rho>0$ there exist $\ve_0=\ve_0(\rho)$ such that for any
$0<\ve<\ve_0$, for any $f\in D(A)$ we have, for some $K>0$}
$$|\E_{\boldq} f (\boldpi(\boldq_{\sm_0}^\ve))-f(\boldpi(\boldq))|<K(\utl(\dt)-\utl(-\dt))^2\ . \eqno(3.28)$$

\

\textbf{Proof.} We have, using Corollary 3.1, that

$$\begin{array}{l}
|\E_{\boldq} f (\boldpi(\boldq_{\sm_0}^\ve))-f(\boldpi(\boldq))|
\\
=\left|\E_{\boldq} f(\tht_{\sm_0}^\ve,
\pi(y_{\sm_0}^\ve))-f(\boldpi(\boldq)) \right|
\\
\play{=\left|\int_0^{2\pi} f(\tht,
\dt)\Prob_{\boldq}(\tht_{\sm_0}^\ve\in d\tht,
y_{\sm_0}^\ve=1+\dt)+\int_0^{2\pi} f(\tht,
-\dt)\Prob_{\boldq}(\tht_{\sm_0}^\ve\in d\tht,
y_{\sm_0}^\ve=-1-\dt)-f(\boldpi(\boldq))\right|}
\\
\play{\leq\left|\dfrac{1}{2\pi}\int_0^{2\pi}
\dfrac{\utl(0)-\utl(-\dt)}{\utl(\dt)-\utl(-\dt)}f(\tht,
\dt)d\tht+\dfrac{1}{2\pi}\int_0^{2\pi}
\dfrac{\utl(\dt)-\utl(0)}{\utl(\dt)-\utl(-\dt)}f(\tht,
-\dt)d\tht-f(\boldpi(\boldq))\right|+K_1 (\utl(\dt)-\utl(-\dt))^2}
\\
\play{=\left|\dfrac{1}{2\pi}\int_0^{2\pi}
\dfrac{\utl(0)-\utl(-\dt)}{\utl(\dt)-\utl(-\dt)}(f(\tht,
\dt)-f(\fo))d\tht\right.-}
\\
\play{ \ \ \ \ \ \ \ \ \ \ \ \left. \dfrac{1}{2\pi}\int_0^{2\pi}
\dfrac{\utl(\dt)-\utl(0)}{\utl(\dt)-\utl(-\dt)}(f(\fo)-f(\tht,-\dt))d\tht
+(f(\fo)-f(\boldpi(\boldq)))\right|+K_1(\utl(\dt)-\utl(-\dt))^2}
\\
\play{\leq
\left|\dfrac{(\utl(0)-\utl(-\dt))(\utl(\dt)-\utl(0))}{\utl(\dt)-\utl(-\dt)}
\left(\dfrac{1}{2\pi}\int_0^{2\pi} \dfrac{f(\tht,
\dt)-f(\fo)}{\utl(\dt)-\utl(0)}d\tht - \dfrac{1}{2\pi}\int_0^{2\pi}
\dfrac{f(\fo)-f(\tht,-\dt)}{\utl(0)-\utl(-\dt)}d\tht\right)\right| +
}
\\
\play{\ \ \ \ \ \ \ \ \ \ \ \ |f(\fo)-f(\boldpi(\boldq))| +
K_1(\utl(\dt)-\utl(-\dt))^2}
\\
\leq K(\utl(\dt)-\utl(-\dt))^2
\end{array}$$
for some $K_1>0$ and $K>0$. We have used the gluing condition (3.10)
and our specified choice of asymptotic order of $\dt$, $\dt'$ and
$\ve$. $\square$

\

\textbf{Lemma 3.10. }\textit{We have, as $\ve,\dt, \dt'$ are small,
for $\boldq_0\in G(\dt)$, that }
$$\left|\E_{\boldq_0}\left[\int_{\sm_0}^{\tau_1}e^{-\lb t}[\lb
f(\boldpi(\boldq_t^\ve))-Af(\boldpi(\boldq_t^\ve))]dt+e^{-\lb
\tau_1}f(\boldpi(\boldq_{\tau_1}^\ve))\right]-f(\boldpi(\boldq_0))\right|\leq
(\utl(\dt)-\utl(-\dt))^2 \ . \eqno(3.29)$$

The proof of this Lemma is essentially the same proof in Lemma 2.6
modified into a two-dimensional version and we omit it.

\

Finally we would like to mention that our boundary condition given
in this section also appears naturally in other model problems. As
an example let consider the following system:

$$\left\{\begin{array}{l}
\play{x_t^\ve=\int_0^t\dfrac{1}{\lb(y_t^\ve)+\ve}dW_t^1 \ ,}
\\
y_t^\ve=|W_t^2| \ .
\end{array}\right. \eqno(3.30)$$

Here $\lb(\bullet)$ is a smooth function on $\R_+$ that vanishes at
$0$ and is strictly positive in $(0,\infty)$; $W_t^1$ and $W_t^2$
are two independent standard Wiener processes on $\R$. Let the
process $z_t^\ve=(x_t^\ve,y_t^\ve)$ on $\R\times \R_+$ be stopped
once it hits the boundary $\{(x,y)\in \R^2: y=R\}$ for some $R>0$.
Let $\tht_t^\ve=x_t^\ve \mod 2\pi$. Let $\boldpi:S^1\times \R_+\ra
\R^2$ be the mapping defined by
$\boldpi(\tht,y)=(y\cos\tht,y\sin\tht)$. For each fixed $\ve>0$, the
process $w_t^\ve=(\tht_t^\ve, y_t^\ve)$ is a diffusion process on
$S^1\times [0,R]$ with normal reflection at the boundary
$\{(\tht,y): y=0\}$ and is stopped once it hits the other boundary
$\{(\tht,y): y=R\}$. Let $m_t^\ve=\boldpi(w_t^\ve)$ (i.e., we glue
all points $\{(\tht,y): y=0\}$). The process $m_t^\ve$ moves within
the disk $B(R)=\{m\in \R^2: |m|_{\R^2}\leq R\}$ and is stopped once
it hits the boundary. In general, this process is \textit{not} a
Markov process. But we expect that, as $\ve \da 0$, this process
$w_t^\ve$ will converge weakly to a Markov process $w_t$ on $B(R)$
with generator $A$ and the domain of definition $D(A)$, defined as
follows: The operator $A$ at points $(\tht,r)$ (we use polar
coordinates, that is, a point $(x,y)\in \R^2$ is represented by
$(r\cos \tht, r\sin \tht)$) with $r\neq 0$ is defined by
$$Af(\tht,r)=\dfrac{1}{2\lb^2(r)}\dfrac{\pt^2}{\pt
\tht^2}f(\tht,r)+\dfrac{1}{2}\dfrac{\pt^2}{\pt r^2}f(\tht,r) \ .
\eqno(3.31)$$ The domain of definition $D(A)$ of the operator $A$
consists of those continuous functions $f$ on $B(R)$ for which
$Af(\tht, r)$ is defined and continuous for $r\neq 0$, the
derivative in $r$ being continuous; such that finite limit
$$\lim\li_{\tht'\ra \tht, r\ra 0+}\dfrac{\pt f}{\pt r}(\tht', r) \eqno(3.32)$$
exists; $$\lim\li_{\tht'\ra \tht, r\ra 0+} Af(\tht', r)
\eqno(3.33)$$ exists and does not depend on $\tht$;
$$\lim\li_{\tht' \ra \tht, r\ra R-}Af(\tht', r)=0 \ ; \eqno(3.34)$$ and
$$\play{\int_0^{2\pi}\lim\li_{\tht'\ra \tht, r\ra 0+}\dfrac{\pt
f}{\pt r}(\tht', r)d\tht}=0 \ . \eqno(3.35)$$

We define, for $f\in D(A)$, $Af(\tht, R)$ as the limit (3.34) and
$Af(O)$ as the limit (3.33).

The weak convergence of $w_t^\ve$ to $w_t$ in
$\contfunc_{[0,T]}(B(R))$ described above shall be a result of fast
motion $x_t^\ve$ running at the local time of the slow motion
$y_t^\ve$ on the boundary $\{(x,y)\in \R\times \R_+: y=0\}$. The
proof of this result shall follow the same method of this section.

\

\section{A conjecture in the general multidimensional case}

In this section we give a conjecture in the general multidimensional
case. Consider the general multidimensional problem (1.8), and for
brevity assume that $\boldb(\bullet)\equiv \mathbf{0}$. That is, the
system has the form

$$\dot{\boldq}^\ve_t=-\dfrac{\grad \lb(\boldq_t^\ve)}{2(\lb(\boldq_t^\ve)+\ve)^3}+
\dfrac{1}{\lb(\boldq_t^\ve)+\ve} \dot{\boldW}_t \ , \
\boldq^\ve_0=\boldq\in \R^d \ . \eqno(4.1)$$

Let us work in a large closed ball $B(R)=\{\boldq\in \R^d:
|\boldq|_{\R^d}\leq R\}$ for some $R>0$, i.e., the process
$\boldq_t^\ve$ is stopped once it hits $\pt B(R)$. Suppose the
friction $\lb(\bullet)$ is smooth and $\lb(\boldq)=0$ for $\boldq$
in some region $G\subset B(R)$ while $\lb(\boldq)>0$ for $\boldq\in
B(R)\backslash [G]$ (here $[G]$ is the closure of $G$ with respect
to the Euclidean metric in $\R^d$). The domain $G\subset B(R)$ is
assumed to be simply connected and to have a smooth boundary $\pt
G$.

Let $\fC$ be a topological space consisting of all points in
$B(R)\backslash [G]$ and one additional point $\fo$. The topology of
$\fC$ contains all the open subsets (in standard Euclidean metric)
in the induced topology of $B(R)\backslash [G]$ and all the open
neighborhoods of $[G]$ in $B(R)$ as the open subsets of $\fC$
containing $\fo$. Let us consider a projection $\boldpi: B(R)\ra
\fC$ defined as follows: for points $\boldq\in B(R)\backslash [G]$
we have $\boldpi(\boldq)=\boldq$ and for points $\boldq\in [G]$ we
have $\boldpi(\boldq)=\fo$. Under the above defined topology for
$\fC$ the mapping $\boldpi$ is continuous. Let
$\boldqtl_t^\ve=\boldpi(\boldq_t^\ve)$ be a stochastic process with
continuous trajectories on $\fC$.

Our conjecture is about the weak convergence, as $\ve\da 0$, of
$\boldqtl_t^\ve$ to some Markov process $\boldqtl_t$ on $\fC$. Below
we give our definition of the latter process but we point out that
we are not clear about the \textit{existence} of it. Our generator
and boundary condition for this process is more or less in the
spirit of martingale problems (see, for example,
\cite[Ch.4]{[Ethier-Kurtz]}). To ensure the \textit{uniqueness} of
solution of martingale problems we need the \textit{existence} of
solution in a nice space of the corresponding PDE with the specified
boundary condition. We are not clear about this yet.

The process $\boldqtl_t$ is identified by its generator $A$ with
domain of definition $D(A)$. For a function $f(\boldqtl)$ on $\fC$
that is continuous on $\fC$ and smooth for $\boldqtl\neq \fo$,
$|\boldqtl|_{\R^d}<R$ we define
$$Af(\boldqtl)=-\dfrac{\grad \lb(\boldqtl)\cdot
\grad f(\boldqtl)}{2\lb^3(\boldqtl)}+\dfrac{1}{2\lb^2(\boldqtl)}\Dt
f (\boldqtl) \ , \eqno(4.2)$$ and at the points $\fo$ and $\boldqtl$
with $|\boldqtl|_{\R^d}$ we define the values of $Af$ as the limits
of the values given by (4.2) (assuming these limits exist). The
domain $D(A)$ is defined as the set of functions $f$ such that
$Af(\boldqtl)=0$ for $|\boldqtl|_{\R^d}=R$, the generalized normal
derivative
$$D_{\utl}f(\boldq)=\lim\li_{\dt \da
0}\dfrac{f(\boldq+\dt \boldn(\boldq))-f(\fo)}{\utl(\boldq+\dt
\boldn(\boldq))} \eqno(4.3)$$ exists for all $\boldq\in \pt G$,
where $\boldn(\boldq)$ is the vector of the outside normal to $\pt
G$, and $\utl(\boldq)$ is some function defined in a neighborhood of
$\pt G$ with $\lim\li_{\boldpi(\boldq)\ra \fo}\utl(\boldq)=0$; and

$$\int_{\pt G}D_{\utl}f(\boldq)d\sm(\boldq)=0 \ . \eqno(4.4)$$
Here $d\sm(\boldq)$ denotes integration with respect to the surface
area on $\pt G$.

\

\textbf{Conjecture.} \textit{The process
$\boldqtl_t^\ve=\boldpi(\boldq_t^\ve)$ converges weakly in the space
$\contfunc_{[0,T]}(\fC)$ as $\ve \da 0$ to the process $\boldqtl_t$
described above.}

\

A further conjecture: we can define the function $\utl$ as
$$\utl(\boldq+\dt\boldn(\boldq))=\int_0^\dt\lb(\boldq+s\boldn(\boldq))ds \eqno(4.5)$$
for $\boldq\in \pt G$ and $\dt>0$ sufficiently small.

\

\textbf{Acknowledgements}: This work is supported in part by NSF
Grants DMS-0803287 and DMS-0854982.

\


\begin{thebibliography}{100}

\bibitem{[Barlow-Pitman-Yor]} Barlow, M., Pitman, J., Yor, M., On
Walsh's Brownian motion, \textit{S\'{e}minaire de Probabiliti\'{e}s
XXIII}, Springer Lecture Notes in Mathematics, \textbf{1372} (1989),
pp.\, 275\,--\,293.

\bibitem{[Dynkin]} Dynkin, E., One dimensional continuous strong Markov processes,
\textit{Theory of Probability and Its Applications} (English
translation), {\bf 4}, 1, 1959, pp.\, 1\,--\,52.

\bibitem{[Ethier-Kurtz]} Ethier, S., Kurtz, T., \textit{Markov processes,
characterization and convergence}, John Wiley and Sons, London,
1986.

\bibitem{[Feller]} Feller, W., Generalized second-order differential
operators and their lateral conditions, \textit{Ill. Journal of
Math.}, {\bf 1} (1957), pp.\, 459\,--\,504.

\bibitem{[F JSP]}  Freidlin, M., Some Remarks on the Smoluchowski-Kramers
Approximation, \textit{Journal of Statistical Physics},
\textbf{117}, No.314, pp.\, 617\,--\,634, 2004.

\bibitem{[FH SK PMA]} Freidlin, M., and Hu, W., Smoluchowski-Kramers
approximation in the case of variable friction, \textit{Journal of
Mathematical Sciences}, \textbf{79}, No.1, November 2011, pp.\,
184\,--\,207, translated from \textit{Problems in Mathematical
Analysis}, \textbf{61}, October 2011 (in Russian).

\bibitem{[FW fish paper]} Freidlin, M., and Wentzell, A., On the Neumann
problem for PDE's with a small parameter and the corresponding
diffusion processes, \textit{Probability Theory and Related Fields},
online, DOI:10.1007/s00440-010-0317-4.

\bibitem{[FW book]} Freidlin, M., and Wentzell, A., \textit{Random
Perturbations of Dynamical Systems}, 2-nd edition, Springer, 1998.

\bibitem{[FW AMS]} Freidlin, M., and Wentzell, A., Random Perturbations of
Hamiltonian Systems, \textit{Mem. of AMS}, {\bf 109} (1994).

\bibitem{[FW Diffusion process on a graph]} Freidlin, M., and Wentzell, A.,
Diffusion processes on graphs and the averaging principle,
\textit{Annals of Probability}, {\bf 21}, 4, 1993, pp.\,
2215\,--\,2245.

\bibitem{[Mandl]} Mandl, P., \textit{Analytical Treatment of One-dimensional
Markov Processes}, Springer, 1968.

\bibitem{[Mochanov]} Mochanov, S., On a problem in the theory of diffusion
processes, \textit{Theory of Probability and Its Applications}
(English translation), {\bf 9}, 1964, pp.\, 472\,--\,477.

\bibitem{[Mochanov-Ostrovskii]} Mochanov, S., and Ostrovskii, E., Symmetric
stable processes as traces of degenerate diffusion processes,
\textit{Theory of Probability and Its Applications} (English
translation), {\bf 14}, 1969, pp.\, 128\,--\,131.

\bibitem{[Revuz-Yor]} Revuz, D., and Yor, M., \textit{Continuous Martingales
and Brownian Motion}, 3.ed., Springer, 1999.

\bibitem{[Ueno I]} Ueno, T., The diffusion satisfying Wentzell's boundary
condition and the Markov process on the boundary I, \textit{Proc.
Japan Acad.}, {\bf 36}, 10, 1960, pp. \, 533\, -- \,538.

\bibitem{[Ueno II]} Ueno, T., The diffusion satisfying Wentzell's boundary
condition and the Markov process on the boundary II, \textit{Proc.
Japan Acad.}, {\bf 36}, 10, 1960, pp. \, 625\, -- \,629.

\bibitem{[Wentzell boundary condition]} Wentzell, A., On boundary
conditions for multidimensional diffusion processes, \textit{Theory
of Probability and Its Applications} (English translation), {\bf 4},
2, 1959, pp.\, 164\,--\,177.

\end{thebibliography}
\end{document}